\documentclass[rep10,leqno,notebook]{article}
\usepackage{amsmath, amssymb,latexsym}

\newcommand{\ds} {\baselineskip=20pt \lineskiplimit=10pt \lineskip=10pt}

\def\polys{\mbox{$\mathbb{F}[k_1,\ldots,k_l,u_1,\ldots,u_{2 n^2 -l} ]$}}

\def\knownpolys{\mbox{$\mathbb{F}[k_1,\ldots,k_l]$}}
\def\I{\mathcal{I}}

\ds

\begin{document}

\newcommand{\tl}[1]{\tilde{#1}}


\newtheorem{defn}{{\bf Definition}}[section]
\newtheorem{thm}{{\bf Theorem}}
\newtheorem{cor}{{\bf Corollary}}
\newtheorem{lem}{{\bf Lemma}}
\newtheorem{pro}{{\bf Proposition}}
\newtheorem{cond}{{\bf Condition}}
\newtheorem{rem}{{\bf Remark}}

\newcommand{\lskip}{\vspace{.18in}}

\newenvironment{proof}{\noindent{\bf Proof:\quad}}{\hfill $\Box$ \lskip}
\newenvironment{preproof}{\noindent{\bf Orientation:\quad}}{\hfill $\Diamond$ \lskip}

\newenvironment{ex}{\lskip \refstepcounter{defn} \noindent{\bf Example
\thedefn}}{\lskip}
\newenvironment{fact}{\lskip \refstepcounter{defn} \noindent{\bf Fact \thedefn}}{}


\def\makeatletter{\catcode`\@=11\relax}


\title{ Using noncommutative
Gr\"obner bases in solving partially prescribed
matrix inverse completion problems }

\author{
 F. Dell Kronewitter\thanks{University of California, San Diego, Math Dept.
 9500 Gilman Dr., San Diego, CA 92093-0112.  Partially supported
by the AFOSR and the NSF.}}

\date{ \today }
\maketitle

\framebox(400,70){\shortstack{
{This paper appears in the following journal: } \\ 
  \large F. Dell Kronewitter, {\bf Using noncommutative
Gr\"obner bases in } \\
\large {\bf solving partially prescribed matrix inverse }
\\ \large {\bf  completion problems, {\it \bf Linear Algebra and its Applications}},
\\ \large
Vol. 338, Pages 171-199, 2001. }}

\begin{abstract}
{
We investigate  the use of noncommutative Gr\"obner bases
in solving partially prescribed matrix inverse completion
problems.     The types of problems considered here are
similar to those in \cite{BLJW}.
There the authors
gave necessary and sufficient conditions for the solution of
a two by two block matrix
completion problem.     Our approach is quite different from
theirs and relies on symbolic computer algebra.

Here we describe a general method by which all block matrix
completion problems of this type may be analyzed if sufficient computational
power is available.   We also demonstrate our method  with an analysis of
all three by three block matrix inverse completion problems
with eleven
blocks
known and seven unknown.   We discover that the solutions
to all such problems are of a relatively simple form.

We then perform a more detailed analysis of a particular
problem from the
31,824  three by three block matrix
completion  problems with eleven
blocks known and seven unknown.   A solution to this problem
of the form derived in \cite{BLJW} is presented.

Not only do we give a proof of our detailed result,
but we describe the strategy
used in discovering our theorem and proof,   since it is somewhat unusual for
these types of problems.
}
\end{abstract}

\newpage

\tableofcontents

\newpage

\section{The problem}
\bigskip

We consider block matrix completion problems similar to that in
~\cite{BLJW}.  Here, we take two partially prescribed, square
matrices, $A$ and $B$, and describe conditions which make it
possible to complete the matrices so that they are inverses of
each other.  That is,  we wish the completed matrices to satisfy
\begin{equation}
\label{inverse-cond}
 A B = I \mbox{ and }  B A = I.
\end{equation}

\subsection{A sample problem }
\label{sec:sample-prob}

An example of such a problem is:  given
matrices $k_1$, $k_2$, $k_3$, and $k_4$,
and
\begin{equation}
\label{sample-2x2}
A =
\left(
\begin{array}{cc}
k_1 & u_2 \\
u_1 & k_2
\end{array}
\right)
\mbox{ and }
B =
\left(
\begin{array}{cc}
u_3 & k_4 \\
k_3 & u_4
\end{array}
\right),
\end{equation}
is it possible to find matrices $u_1$, $u_2$, $u_3$, and $u_4$
such that equation (\ref{inverse-cond}) is satisfied? The answer
to this question, due to \cite{BLJW},
 is given in Section \ref{BJLW-soln}.

We now describe our problem in detail.

\subsection{The general block matrix inverse completion problem }
\label{sec:general-mtx-complete}

We
begin by partitioning two matrices,  $A$ and $B$, whose entries
are elements in an arbitrary field $\mathbb{F}$,  conformally for
matrix multiplication into $n$ by $n$ block matrices.    Next,
 we choose $l$ of these
blocks to be known and $2 n^2 - l$ to be unknown.   We give some
conditions on the known matrices, which may be expressed
algebraically,  such as invertibility or self-adjointness. We will
now define our problem.

    We ask if it is possible to fill in the $2 n^2 - l$ unknown
blocks so that equation (\ref{inverse-cond}) is satisfied and seek
to derive formulas for these matrices in terms of the prescribed
blocks.  To be more specific, we might even call this problem the
{\em purely algebraic partially prescribed
matrix inverse completion
problem}.
The solution to such a problem will be a set of matrix equations
in the known and unknown submatrices.

\section{The solution } \label{tri-solution}

In general,  it is not known how to solve a system of matrix
equations where several of the matrices are unknown. Unknown
matrices can appear in matrix equations in several ways, of which
some are more computationally acceptable than others.   We will
analyze these forms and classify  certain solution sets. In this
section,  as well as throughout the paper,
noncommutative variables $k_i$
will be considered initially known, and noncommutative variables
$u_i$ will be considered initially unknown.

We first recall some standard definitions.
If $\mathbb{F}[k_1,\ldots,k_l,u_1,\ldots,u_{2 n^2 -l}]$ is a
noncommutative, free algebra over the field $\mathbb{F}$ of characteristic
$0$,  then a subset $\mathcal{I}$ will be called an {\em \bf ideal} if
$f,g \in \I$ implies that $l_f f r_f + l_g g r_g \in \I$ for $l_f,r_f,l_g$, and
$r_g$ arbitrary elements
of  $\mathbb{F}[k_1,\ldots,k_l,u_1,\ldots,u_{2 n^2 -l}]$.
Given a finite set of polynomials $\{ p_1,\ldots,p_m \}$ in
$\mathbb{F}[k_1,\ldots,k_l,u_1,\ldots,u_{2 n^2 -l}]$, we say that $\I$
{\em \bf is generated by} $\{ p_1,\ldots,p_m \}$ if $\I$ is the smallest ideal
containing $\{ p_1,\ldots,p_m \}$.

\subsection{A good,  triangular solution } \label{sec:goodanswer}

The following definition is useful in identifying a class of
problems, particularly those consisting of matrix equations,
which are usually computationally tractable. This will be
demonstrated in Section~\ref{sec:solvegood}.  It is, in general,
impossible to verify condition (\ref{non-deg-cond}) in
Definition~\ref{f-back} below. For this reason,  we give two
versions of our definition:  a weaker, computable, version and a
stronger, in general incomputable, version. We will later
introduce approximations to condition (\ref{non-deg-cond}) which
can be verified.
\begin{defn}
\label{f-back} Let ${\cal I}$ be an ideal in $\polys$. We say that
${\cal I}$ {\em \bf can be made weakly formally backsolvable} if
there exists a bijective map $$\sigma : \{1,\ldots,2 n^2 - l\}
\rightarrow \{ 1,\ldots,2 n^2 - l \}$$
  and a
finite set $G$ of polynomials which generates ${\cal I}$ such that
\begin{eqnarray}
& G = G_0 \cup G_1 \cup \ldots \cup G_{2 n^2 - l}
& \nonumber
\end{eqnarray}
where $G_0 = G \cap \knownpolys$
    and
$$G_i \subset \mathbb{F}[k_1,\ldots,k_l,u_{\sigma(1)},\ldots,u_{\sigma(i)}]
    \setminus \mathbb{F}[k_1,\ldots,k_l,u_{\sigma(1)},\ldots,u_{\sigma(i-1)}],
\qquad G_i \neq \emptyset,$$
for
    $1 \leq i \leq 2 n^2 - l$.  If, in addition, we have the
    condition,
\begin{equation}
\mbox{no proper subset of $G$ generates ${\cal I}$},
\label{non-deg-cond}
\end{equation}
we say $\mathcal{I}$ {\em \bf can be made formally backsolvable}.
We say that the set $G$ {\em \bf is formally backsolvable.}
\end{defn}

Definition~\ref{f-back}, for many people, will be more intuitive
in an expanded notation. What follows is an intuitive and expanded
notation for weakly formally backsolvable. Indeed, the set $G$ of
polynomials $G = \{q_{0,1},\ldots,q_{2 n^2 - l, m_{2 n^2 -l }}\}$
has the form
\begin{eqnarray}
q_{0,1}(k_1,\ldots,k_l) &=&0 \label{tri-1} \\
q_{0,2}(k_1,\ldots,k_l) &=&0
\\
 & \vdots& \nonumber \\
q_{0,m_0}(k_1,\ldots,k_l) &=&0 \label{tri-known}  \\
 & & \nonumber \\
 q_{1,1}(k_1,\ldots,k_l,{\bf u_{\sigma(1)}})&=&0  \label{def1:u1}
 \\
 & \vdots & \nonumber \\
 q_{1,m_1}(k_1,\ldots,k_l,{\bf u_{\sigma(1)}})&=&0 \label{defn:u1} \\
 q_{2,1}(k_1,\ldots,k_l,u_{\sigma(1)},
{\bf u_{\sigma(2)} }) &=& 0 \label{def1:u2} \\
& \vdots &  \nonumber \\
q_{2,m_2}(k_1,\ldots,k_l,u_{\sigma(1)},
{\bf u_{\sigma(2)} })&=&0  \label{defn:u2}  \\
q_{3,1}(k_1,\ldots,k_l,u_{\sigma(1)},{u_{\sigma(2) },
{\bf u_{\sigma(3)} }   })&=&0 \nonumber \\
 & \vdots& \nonumber \\
 q_{2 n^2 - l - 1, m_{2 n^2 - l - 1} }
(k_1,\ldots,k_l,u_{\sigma(1)},u_{\sigma(2)},\ldots,{\bf
u_{\sigma(2n^2-l-1)}})&=&0 \nonumber \\
 q_{2 n^2 - l, 1 }
(k_1,\ldots,k_l,u_{\sigma(1)},u_{\sigma(2)},u_{\sigma(3)},u_{\sigma(4)},\ldots,{\bf
u_{\sigma(2n^2-l)}})&=&0 \nonumber \\
& \vdots & \nonumber \\
 q_{2 n^2 - l, m_{2n^2 -l} }
(k_1,\ldots,k_l,u_{\sigma(1)},u_{\sigma(2)},u_{\sigma(3)},u_{\sigma(4)},\ldots,{\bf
u_{\sigma(2n^2-l)}})&=&0 \label{tri-last} \,
\end{eqnarray}
where the $k_i$ are known, $m_i > 0$, the $u_{\sigma(i)}$ are
unknown, and $\sigma$ is a permutation map on integers $1$ to $2
n^2 -l$. The first subscript of the polynomials $q_{r,m}$
indicates the number of unknowns allowed in the equation, as well
as implying the existence of the bolded unknown ${\bf
u_{\sigma(r)}}$.

We refer to equations (\ref{tri-1}-\ref{tri-known})
which contain only knowns
as {\bf compatibility conditions on the knowns}.
These equations, in only the known variables, must hold
 if a completion is possible.
Equations (\ref{def1:u1}-\ref{tri-last}) containing unknowns we
call {\bf equations triangular in the unknowns} due to the
triangular structure exhibited by the unknown variables.

Usefulness for computation is discussed
in Section~\ref{sec:solvebetter}.

\subsection{Decoupled solutions }
\label{sec:betteranswer}

The following definitions are useful in identifying two classes of
problems, particularly those consisting of matrix equations,
which are usually even more computationally tractable than the
formally backsolvable form,  introduced in the previous section.
We call the two classes essentially decoupled and formally decoupled.
It
is not always possible to verify condition
(\ref{non-deg-cond-better}) in Definition~\ref{def:ess-decoup}
below.
For this reason,  we give two versions of each class definition:
a weaker, computable, version and a stronger, in general
incomputable, version. Later, we will introduce approximations to
condition (\ref{non-deg-cond-better}), which can be verified.
\begin{defn}
\label{def:ess-decoup}
    Let ${\cal I}$ be an ideal in $\polys$.
Let there exist $j$, $1\leq j\leq 2 n^2 -l$;
    an injective map $\sigma : \{1,\ldots,j\} \rightarrow
    \{ 1,\ldots,2 n^2 -l  \}$; a bijective map $\tau:
    \{ j+1,\ldots,2 n^2 - l\} \rightarrow
    \{1,\ldots,2 n^2 -l \} \setminus \mbox{image}(\sigma) $; a set $G^*$
    of polynomials;
  and a
finite  set $G$ of polynomials which generates ${\cal I}$ such
that
\begin{eqnarray}
& G = G_0 \cup G_1 \cup \ldots \cup G_j \cup
    \{ u_{\tau(i)} - g_i(k_1,\ldots,k_l,u_{\sigma(1)},\ldots,u_{\sigma(j)})
    : j+1 \leq i \leq 2 n^2 -l \} & \nonumber \\
& \cup \: G^* & \nonumber
\end{eqnarray}
where $G_0 = G \cap \knownpolys$
    and
$$G_i \subset \mathbb{F}[k_1,\ldots,k_l,u_{\sigma(i)}]
    \setminus \mathbb{F}[k_1,\ldots,k_l], \qquad G_i \neq \emptyset,$$
for
    $1 \leq i \leq j$.

We now make the following definitions concerning the ideal introduced above.
\begin{enumerate}
\item{\em We say that the polynomial ideal $\mathcal{I}$
{\bf can be
weakly essentially decoupled}
and that the set of polynomials $G$ {\em \bf is weakly essentially
decoupled.}
}

\item{\em If $G^* = \emptyset$, then we say that the polynomial
ideal $\mathcal{I}$ {\em \bf can be
weakly formally decoupled}  and that the set of
polynomials $G$ {\bf is weakly formally
decoupled. }
}
\bigskip

\noindent
An important nondegeneracy condition is
\begin{equation}
\label{non-deg-cond-better}
 \mbox{no proper subset of $G$
generates ${\cal I}$.}
\end{equation}

\item{\em If condition (\ref{non-deg-cond-better}) holds,  then we say
 that the polynomial
ideal $\mathcal{I}$ {\bf can be
essentially decoupled} and
 that the
set of polynomials $G$ {\bf is
essentially decoupled}.}

\item{\em If $G^* = \emptyset$ and condition
(\ref{non-deg-cond-better}) holds,  then we say
 that the polynomial
ideal $\mathcal{I}$ {\bf can be
formally decoupled}
and
 that the
set of polynomials $G$ {\bf is
formally decoupled}.
}
\end{enumerate}
\end{defn}

\medskip

Definition~\ref{def:ess-decoup}, for many people, will be more
intuitive in an expanded notation. What follows is an intuitive
and expanded notation for weakly decoupled.  Indeed, the set $G$
of polynomials,
$$G = \{q_{0,1},q_{0,2},\ldots,q_{j,m_j},q_{s_j+1}-u_{\tau(j+1)},
\ldots,q_{s_j+2 n^2-l -j}-u_{\tau(2 n^2 -l )}\}$$ in the formally
decoupled case or
$$G = \{ q_{0,1},q_{0,2},\ldots,q_{j,m_j},q_{s_j+1}-u_{\tau(j+1)},
\ldots,q_{s_j+2 n^2-l -j}-u_{\tau( 2 n^2 -l )}, q_{ s_j + 2 n^2 -l
-j + 1 }, \ldots, q_{s_0} \}$$
 in the essentially decoupled case,
has the form
\begin{eqnarray}
q_{0,1}(k_1,\ldots,k_l) &=&0 \label{flat-1} \\
q_{0,2}(k_1,\ldots,k_l)
&=&0 \\
 & \vdots& \nonumber \\
q_{0,m_0}(k_1,\ldots,k_l) &=&0 \label{endknowns} \\
& & \nonumber \\
q_{1,1}(k_1,\ldots,k_l,{\bf u_{\sigma(1)}})&=&0 \label{beginunk}
\\
 & \vdots& \nonumber \\
q_{1,m_1}(k_1,\ldots,k_l,{\bf u_{\sigma(1)}})&=&0
\\
q_{2,1}(k_1,\ldots,k_l,{\bf u_{\sigma(2)}})&=&0
\\
 & \vdots& \nonumber \\
q_{j-1,m_{j-1}}(k_1,\ldots,k_l,{\bf u_{\sigma(j-1)}})&=&0  \\
q_{j,1}(k_1,\ldots,k_l,{\bf u_{\sigma(j)}})&=&0  \\
& \vdots & \nonumber \\
q_{j,m_j}(k_1,\ldots,k_l,{\bf u_{\sigma(j)}})&=&0 \label{end-det} \\
& & \nonumber \\
{\bf u_{\tau(j+1)}}&=&q_{s_j+1}(k_1,\ldots,k_l,
u_{\sigma(1)},\ldots,u_{\sigma(j)}) \label{single} \\
 & \vdots& \nonumber \\
{\bf u_{\tau(2 n^2-l)}}&=&q_{s_j + 2 n^2 -l -j }(k_1,\ldots,k_l,
u_{\sigma(1)},\ldots,u_{\sigma(j)}) \label{end-singles}
\end{eqnarray}
\begin{quotation}
{\em
\noindent
The following equations, which we call
{\bf compatibility conditions on the unknowns}, will not
exist in the formally decoupled or weakly formally decoupled cases,
but might occur in the essentially decoupled case.}
\end{quotation}

\vspace{-.5cm}

\begin{eqnarray}
q_{s_j + 2 n^2 -l -j + 1 }(k_1,\ldots,k_l,u_1,\ldots,u_{2 n^2 -l})
& = & 0
\label{compat-full} \\
 & \vdots& \nonumber \\
q_{s_0}(k_1,\ldots,k_l,u_1,\ldots,u_{2 n^2 - l}) & = & 0
\label{flat-last}
\end{eqnarray}

\noindent
 In the above system of equations $s_j = \sum_{i=0}^j m_i$ with $m_i > 0$
  and
 $s_0 = \mid \! G \!  \mid$.
The first non-zero subscript of the doubly subscripted polynomials
$q_{r,m}$ indicates the occurrence of one and only one unknown
$u_{\sigma(r)}$.

 Equations of the form
(\ref{single}-\ref{end-singles}) will be referred to as {\bf
singletons}.   A singleton equation is characterized by the fact
that there is a single instance of an unknown variable which does
not occur in equations (\ref{beginunk}-\ref{end-det}).   This
unknown variable appears in the singleton equation as a monomial
consisting of only itself. The {\bf singleton variable} is the
left hand side of equations (\ref{single}-\ref{end-singles}),
$u_{\tau(j+1)},\ldots,u_{\tau(2 n^2 -l)}$. Singleton equations in
Definition~\ref{def:ess-decoup} are $u_{\tau(i)} -
g_i(k_1,\ldots,k_l,u_{\sigma(1)},\ldots,u_{\sigma(j)}) = 0$.

The singleton equation has a very attractive form for a human who
wishes to find polynomials in few unknowns.   Given an equation in
knowns and unknowns ${\cal E}$, it allows one to eliminate the
unknown singleton variable,  for example $u_{\tau(j+1)}$, from
${\cal E}$ by replacing instances of the unknown indeterminate
with its equivalent polynomial representation, in the example case
$q_{s_j+1}$. After this substitution has been performed, the
equation  ${\cal E}$ will not contain the unknown singleton
variable.

As in the formally backsolvable case,  we have {\bf compatibility
conditions on the knowns} (\ref{flat-1}-\ref{endknowns}).
These equations, in only the known variables, must hold
 if a completion is possible.
All unknown variables $u_1,\ldots,u_{2 n^2 -l}$ in equations
(\ref{flat-1}-\ref{flat-last}), which are not singleton unknowns,
appear in equations (\ref{beginunk}-\ref{end-det}) without any
other unknown variables.   Therefore, we think of this
system of equations as
decoupled.    We call equations (\ref{beginunk}-\ref{end-det})
{\bf equations in one unknown}.

In the formally decoupled case,  the coupling
compatibility conditions
on the unknowns (\ref{compat-full}-\ref{flat-last}) are absent.
 This is obviously a better form of solution than the
essentially decoupled form,
since any solutions for $u_{\sigma(1)},\ldots,u_{\sigma(j)}$ will
do. In the essentially decoupled case, one must verify potential solutions
for $u_{\sigma(1)},\ldots,u_{\sigma(j)}$ with
equations (\ref{compat-full}-\ref{flat-last}).

Notice that these decoupled solution forms,
essentially decoupled and formally decoupled,
satisfy the formally backsolvable criteria.   We have the following
set inclusion relationship:
\begin{equation*}
\mbox{Formally decoupled} \subseteq
\mbox{Essentially decoupled} \subseteq
\mbox{Formally backsolvable }
\end{equation*}

\subsection{A sample answer } \label{BJLW-soln}

Here, we give the solution to the problem presented in
Section~\ref{sec:sample-prob}. In \cite{BLJW},  it was shown that,
for invertible $k_i$,  the matrices $A$ and $B$ defined in
(\ref{sample-2x2}) satisfy (\ref{inverse-cond}) if and only if
the unknown submatrix $u_4$ satisfies the following relation
\begin{equation}
\label{def-sample}
 u_4 k_2 u_4 = u_4 + k_3 k_1 k_4.
\end{equation}
The other unknown submatrices are then given in terms
of $u_4$:
\begin{eqnarray}
u_1 &=& k_4^{-1} - k_2 u_4 k_4^{-1} \label{u1} \\ u_2 &=& k_3^{-1} -
k_3^{-1} u_4 k_2  \\ u_3 &=& k_4 k_3 u_4 k_4^{-1} k_1^{-1}  \label{u3}
\end{eqnarray}
This answer contains no compatibility conditions on the knowns.
Equation (\ref{def-sample}) is an equation in one unknown $u_4$.
The remaining equations (\ref{u1}-\ref{u3}) are singletons.
Therefore, the equations associated with this
matrix completion problem can be formally decoupled.
In such circumstances, we would say that the problem can
be formally decoupled.
  In the
language of Definition~\ref{def:ess-decoup}, $G_0 = \emptyset$,
$G^* = \emptyset$,
$j=1$,
$\sigma(1) = 1$, $G_1 = \{ (\ref{def-sample}) \}$, $\tau(2) = 2$,
$\tau(3) = 3$, and $\tau(4) = 4$.

This main theorem of \cite{BLJW} is simpler,  from a
computational perspective, than the results
we are presenting here,
since (\ref{sample-2x2}) contains fewer equations and
fewer variables. In addition to being proven via traditional
methods, the main theorem in \cite{BLJW} was also proven using
noncommutative Gr\"obner methods in \cite{HS}.

\section{Main results on 3x3 matrix inverse completion problems }
\label{sec:mainresults}

We have performed extensive analysis of the 3x3 block matrix
inverse completion problem.
 In particular, we have concentrated on  the problem described in
Section~\ref{sec:general-mtx-complete} where $n$ is three and $l$
is eleven.
We have assumed in our detailed analysis that all
eleven known blocks are invertible.

We begin by noticing that if one matrix problem is a permutation
of another,  then a solution to one transforms to a solution to the other.
We then define a property that characterizes
certain matrix completion problems
which we call {\em strongly undetermined}.
We will present a classification result which characterizes the solutions
to our seven unknown, eleven known matrix completion problems.
In this result, strongly undetermined problems are the worst behaved class.

This section also includes a detailed result on a particular
matrix completion problem,  which is in the same spirit as the \cite{BLJW}
result described in
Section \ref{BJLW-soln}.
We will conclude this section by
outlining  a method by which
one may find a numerical solution to a matrix completion problem
using our symbolic solution.  Much of the work in this
paper appears in the Ph. D. thesis of Dell Kronewitter, \cite{DellThesis}.

\subsection{Configurations and permutations }
\label{configsNperms}

In our investigations of 3x3 block matrix completion problems,
we will refer to a {\bf configuration} as a classification of
blocks  into knowns and unknowns.
We will specify a configuration with $k$'s and $u$'s.
For example,
\begin{equation}
 A =
\left(
\begin{array}{ccc}
 u_1 & k_1 & u_2 \\
k_2 & k_3 & u_3 \\
k_4 & u_4  & k_5
\end{array}
\right) \mbox{ and }
B =
 \left(
\begin{array}{ccc}
k_6 & k_7 & u_5 \\
k_8 & k_9  & u_6 \\
k_{10} & k_{11} & u_7
\end{array}
\right)
\end{equation}
is a configuration.

A 3x3 {\bf permutation matrix} is a 3x3 matrix consisting of three  $1$'s
and six $0$'s.   No two $1$'s may appear in the same row or the same column.
There are, of  course, six such matrices.

For a given 3x3 configuration of knowns and unknowns,  one may
apply 3x3 (block) permutation matrices, $\Pi$ and $\Psi$,  to $A$
and $B$  to get $\Pi^{-1} A \Psi$ and $\Psi^{-1} B \Pi$  and
obtain at most 36 other equivalent configurations. That is,
$$A B = I \mbox{ and }  B A = I$$
if and only if
$$ \Pi^{-1} A
\Psi \:
\Psi^{-1} B \Pi = I \mbox{ and }
\Psi^{-1} B \Pi \:  \Pi^{-1} A
\Psi = I.$$
In describing solutions $A$ and $B$ to this problem,
we will only give one member
of a particular equivalence class $\{ \Pi^{-1} A \Psi,  \Psi^{-1}
B \Pi\}$.

\subsection{Strongly undetermined }

Assume that the pair of block matrices, $A$ and $B$,  are
partitioned into known and unknown blocks that are  compatible for
matrix multiplication.

\begin{defn}
$A$ and $B$ are said to be {\em \bf strongly undetermined} if
there exists an entry of the block matrices
$A B$ or $B A$, which is a polynomial
consisting entirely of unknown blocks.
\end{defn}

Notice that $A$ and $B$ being strongly undetermined is equivalent
to the existence of {\em both} an entire row (column) of unknown
blocks in $A$ {\em and} an entire column (row) of unknown blocks
in $B$.  For example,  the following configuration of known and
unknown blocks is strongly undetermined.
\begin{equation}
\label{strng-und}
 A =
\left(
\begin{array}{ccc}
{\bf u_1} &{\bf  u_2} & {\bf u_3} \\
k_1 & k_2 & k_3 \\
k_4 & k_5 & k_6
\end{array}
\right) \mbox{ and }
B =
 \left(
\begin{array}{ccc}
k_7 & k_8 & {\bf u_4} \\
k_9 & u_5  & {\bf u_6} \\
k_{10} & k_{11} & {\bf u_7 }
\end{array}
\right)
\end{equation}
The product of these two matrices,  $A B$, has the following form.
$$
 \left(
\begin{array}{ccc}
u_1 k_7 + u_2 k_9 + u_3 k_{10} &
u_1 k_8 + u_2 u_5 + u_3 k_{11} &
{\bf u_1 u_4 + u_2 u_6 + u_3 u_7 } \\
k_1 k_7 + k_2 k_9 + k_3 k_{10} &
k_1 k_8 + k_2 u_5 + k_3 k_{11}  &
k_1 u_4 + k_2 u_6 + k_3 u_7  \\
k_4 k_7 + k_5 k_9 + k_6 k_{10} &
k_4 k_8 + k_5 u_5 + k_6 k_{11} &
k_4 u_4 + k_5 u_6 + k_6 u_7
\end{array}
\right)
$$
Since the upper right entry (in boldface) is
a polynomial made up entirely of unknown blocks,
 configuration (\ref{strng-und}) is strongly
undetermined.

\subsection{A class of 31,824 3x3 matrix inverse completion problems }

In our investigations,  we have analyzed (via computer)
a certain collection of 3x3
matrix completion problems.
Two 3x3 block matrices have a total of 18 entries. We have
analyzed those
 which have seven unknown and 11 known
blocks and do not have the strongly undetermined property. We have
chosen to put efforts into this ratio of known to unknown blocks
because we believe Theorem~\ref{particular-thm},  the initial
subject of our research,  to be surprising,  and yet lack the
computational resources to study all 3x3 matrix completion
problems,  or even one 4x4 matrix completion problem.
Section~\ref{sec-recipe} describes how the motivated researcher
with unlimited computational power can go about analyzing a
block matrix problem of any size  of the type addressed in this
paper.

The following theorem shows that all of our seven unknown, 11 known
block matrix completion problems (which are not strongly
undetermined) have  particularly nice solutions.

\begin{thm}
\label{big-thm} Let $A$ and $B$ be three by three block matrices
such that 11 of the 18 blocks are known and seven are unknown. Let
the known blocks be invertible.   The corresponding
partially prescribed inverse
completion problems may be classified as follows.
\begin{enumerate}
\item{
If the configuration of unknown blocks is not strongly
undetermined, and is not of the form given in (\ref{bad}) or a
permutation of such configuration, then the partially prescribed
inverse completion problem {\em can be weakly essentially
decoupled}, in the sense of Definition~\ref{def:ess-decoup}. }
\item{
Problem (\ref{bad}) {\em can be made weakly
formally backsolvable}  in the
sense
 of Definition~\ref{f-back}.
 Thus,  all but the strongly undetermined cases {\em can be made
weakly
formally backsolvable.
}}
\end{enumerate}
 These
answers,
that is the resulting weakly formally backsolvable or weakly essentially
decoupled systems of equations,
 satisfy a  technical non-redundancy condition, {\em
compatibility
3-nondegeneracy}, which will be defined in
Section~\ref{sec:nondeg} once we have built up our Gr\"obner
machinery.
\end{thm}

\medskip

The exceptional case mentioned in Theorem~\ref{big-thm} is
\begin{equation}
\label{bad}
 A =
\left(
\begin{array}{ccc}
k_1 & k_2 & k_3 \\
k_4 & k_5 & k_6 \\
u_1 & u_2 & u_3
\end{array}
\right) \mbox{ and }
B =
 \left(
\begin{array}{ccc}
k_7 & k_8 & u_4 \\
k_9 & u_5  & k_{10} \\
k_{11} & u_{6} & u_7
\end{array}
\right)
\end{equation}

The proof of this theorem,  which requires noncommutative symbolic
software, will be given in Section~\ref{big-thm-proof}.
Answers to the individual problems,   which consist of sets of
polynomials similar to that found in equations
(\ref{def-sample}-\ref{u3}),  can be found on the internet at \\
\verb+http://arXiv.org/abs/math.LA/0101245+.

\subsection{Detailed analysis of a particular 3x3 matrix inverse
completion
problem }
\label{sec:part-thm}

We now give a closer analysis than that given in the
last section  of a particular matrix inverse completion
problem,  which satisfies the assumptions of Theorem~\ref{big-thm}.
We show how someone, interested in a particular matrix
completion problem,  might arrive at a  finer analysis of the problem
instead of the rather terse conclusion given in Theorem~\ref{big-thm}.
Our goal in this section is to present a short,
computationally simple set of formulas which give the solution to
a particular partially prescribed inverse matrix completion problem.
Our conclusions will have the same flavor as those presented in
\cite{BLJW}.

We will analyze a particular problem from those addressed in
Theorem~\ref{big-thm},   the known/unknown configuration:
\begin{equation}
\label{central-prob}
A =
\left(
\begin{array}{ccc}
a & {\underline t} & b \\
{\underline u} & c & {\underline v} \\
d & {\underline w} & e
\end{array}
\right)
\mbox{ and }
B =
 \left(
\begin{array}{ccc}
{\underline x} & f & g \\
h & {\underline y} & i \\
j & k & {\underline z }
\end{array}
\right)
\end{equation}
or the equivalent permuted form,
\begin{equation*}
A =
 \left(
\begin{array}{ccc}
a & b & {\underline t}  \\
d & e & {\underline w} \\
{\underline u} & {\underline v} & c
\end{array}
\right)
\mbox{ and }
B =
 \left(
\begin{array}{ccc}
{\underline x} & g & f \\
j & {\underline z} & k \\
h & i & {\underline y }
\end{array}
\right),
\end{equation*}
where $a$ through $k$ are known and invertible block matrices,
and the underlined ${\underline t}$ through ${\underline z}$
are unknown block matrices.

\begin{thm}
\label{particular-thm}
Given $A$ and $B$ as in (\ref{central-prob}) with invertible knowns
$a$,$b$,$c$,$d$,$e$,$f$,$g$,$h$,$i$,$j$, and $k$, as well as the invertibility
of the matrix  made up of the outer known blocks of $A$ in
(\ref{central-prob}),
\begin{equation}
\left(
\begin{array}{cc}
a & b \\
d & e
\end{array}
\right),    \label{outer}
\end{equation}
 then $A B = I$ and $B A = I$ if and only if the knowns satisfy the
following compatibility conditions:
\begin{eqnarray}
\tl{p} (d a^{-1} - e b^{-1}) & =&
(a^{-1} b - d^{-1} e)  \tl{q}
\label{cc-1}  \\
(d a^{-1} - e b^{-1})^{-1}
\tl{q}& = & (d a^{-1} - e b^{-1})^{-1}
\tl{q} e
(d a^{-1} -
 e b^{-1})^{-1}  \tl{q}  \label{cc-2} \\
& & + (d a^{-1} - e b^{-1})^{-1}
\tl{q}
d g
+ j b
(d a^{-1} - e b^{-1})^{-1}  \tl{q}
-k c i + j a g  \nonumber  \\
 \tl{p}
(a^{-1} b - d^{-1} e)^{-1} & =&
 \tl{p}
(a^{-1} b - d^{-1} e)^{-1}
e
\tl{p}
(a^{-1} b - d^{-1} e)^{-1}   \label{cc-3}  \\
& &
+  \tl{p}
(a^{-1} b - d^{-1} e)^{-1} d g
 + j b  \tl{p}
(a^{-1} b - d^{-1} e)^{-1}
- k c i + j a g   \nonumber
\end{eqnarray}
where
\begin{eqnarray}
\tl{p} & \triangleq & (- a ^{-1} h^{-1} i + a^{-1} b j h^{-1} i -
d^{-1} - d^{-1} e j h^{-1})  \label{def} \\
\tl{q} & \triangleq & ( - k f^{-1} a^{-1}
+ k f^{-1} g d a^{-1}  -  b^{-1} - k f^{-1} g e b^{-1} )
\end{eqnarray}
The unknown matrices can then be given as:
\begin{equation}
 z =
(a^{-1} b - d^{-1} e)^{-1}
(- a ^{-1} h^{-1} i + a^{-1} b j h^{-1} i - d^{-1} - d^{-1} e j h^{-1})
=  (a^{-1} b - d^{-1} e)^{-1} \tl{p}  \label{z-rel2}
\end{equation}
or equivalently
\begin{equation}
 z =
( - k f^{-1} a^{-1}
+ k f^{-1} g d a^{-1}  -  b^{-1} - k f^{-1} g e b^{-1} )
(d a^{-1} - e b^{-1})^{-1}
= \tl{q}  (d a^{-1} - e b^{-1})^{-1}  \label{z-rel1}
\end{equation}
and then
\begin{eqnarray}
 t &=& - a g i^{-1} - b z i^{-1}  \label{t-rel}\\
 u &=& - k^{-1} j a - k^{-1} z d \\
 v &=& k^{-1} - k^{-1} j b - k^{-1} z e \\
 w &=& i^{-1} - d g i^{-1} - e z i^{-1} \\
 x &=& a^{-1} + f k^{-1} j - g d a^{-1} + f k^{-1} z d a^{-1}\\
 y &=& c^{-1} k^{-1} j a f + c^{-1} k^{-1} j b k +
c^{-1} k^{-1} z d f + c^{-1} k^{-1} z e k  \label{y-rel}
\end{eqnarray}
\end{thm}

\bigskip

\noindent

This answer consists of three
 compatibility conditions on the knowns (\ref{cc-1}-\ref{cc-3}),
 an equation in one unknown (\ref{z-rel2}) or (\ref{z-rel1}), and
 singletons (\ref{t-rel}-\ref{y-rel}), and is, therefore,
 {\em formally decoupled}.

The proof of this theorem will be given in
Section~\ref{sec:part-proof}.  We mention some key points
about how the proofs in Theorem~\ref{big-thm}
and Theorem~\ref{particular-thm} compare in
Section~\ref{sec:strength}.
Solutions to all of the problems (configurations) addressed
in Theorem~\ref{big-thm} can be found via the
internet at \\
\verb+http://arXiv.org/abs/math.LA/0101245+.
These solutions consist of a formatted list of equations
in both \LaTeX\  and {\em Mathematica} form.

\subsection{Numerical solutions to matrix equations}
\label{num-solve}

Now we want to see how symbolic solutions
to matrix completion problems  may
be applied to numerical problems in order to
find numerical
matrix completions.
Let us see what is involved in the numerical
solution
of a matrix completion problem with
configuration (\ref{central-prob})
which was solved symbolically
in Theorem 2.   Assume, for
example, that the matrices $A$ and
$B$ are 12x12 and therefore each block is a 4x4 matrix.   That is,
we are given 4x4 matrices $a,b,c,d,e,f,g,h,i,j,$ and $k$
and are trying to find 4x4 matrices
$t,u,v,w,x,y,$ and $z$ which will form the completed inverse.
We now apply Theorem 2 and the formally decoupled set of
equations (\ref{cc-1}-\ref{y-rel}).

 The first step is
to determine whether the matrices $a,b,c,d,e,f,g,h,i,j,$ and $k$
are compatible.
That is,  they must satisfy equations (\ref{cc-1}-\ref{cc-3}).
If these equations are satisfied,  then
the next step is to determine a value for our unknown $z$ using
either equation (\ref{z-rel2}) or (\ref{z-rel1}).
Since we have assumed the invertibility of $a^{-1} b - d^{-1} e$,
the coefficient
of $z$,  this step cannot fail.
 Finally,  one can determine values for the 4x4 matrices
$t,u,v,w,x,$ and $y$
 from the singleton equations (\ref{t-rel}-\ref{y-rel}).
If all of these steps have occurred successfully, then we have
formed our inverse matrix completion.
This illustrates very general behavior which we now describe.

\subsubsection{Solving a decoupled system of matrix equations}
\label{sec:solvebetter}

 A general decoupled system of matrix equations is made up of {\em
compatibility conditions on the knowns}, equations
(\ref{flat-1}-\ref{endknowns});  {\em equations in one unknown},
equations (\ref{beginunk}-\ref{end-det}); {\em singletons},
equations (\ref{single}-\ref{end-singles});  and possibly {\em compatibility
conditions on the unknowns}, equations
(\ref{compat-full}-\ref{flat-last}).

Given such a decoupled set of matrix equations,
one can first verify that a completion may
be possible by verifying equations (\ref{flat-1}-\ref{endknowns})
containing only the given (known) matrices.
Then
 one can use equations
(\ref{beginunk}-\ref{end-det})
$$q_{1,1}=0,\ldots,q_{j,m_j}=0$$ to {\em simultaneously} solve
for the (possibly non-unique) matrices ${\bf
u_{\sigma(1)},\ldots,u_{\sigma(j)}}$ or to determine that
solutions do not exist. Notice that this may constitute a
difficult numerical problem by itself, especially if the matrices
under consideration are of large dimension. It is then a simple
matter to find matrices ${\bf u_{\tau(j+1)},\ldots,u_{\tau(2
n^2-l)}}$ by evaluating polynomials
$$q_{s_j+1},\ldots,q_{s_j + 2 n^2 -l -j }.$$ The boldface $u$'s
in equations (\ref{beginunk}-\ref{end-det}) indicate the unknown
matrix being solved for in each step.

In the essentially decoupled case,
 one must
check  that
the solutions $u_1,\ldots,u_{2n^2-l}$ which this procedure
derives
are acceptable by validating compatibility equations
(\ref{compat-full}-\ref{flat-last}).
Of course, an advantage of decoupled equations
is that their solution may easily be parallelized.

\subsubsection{Solving a formally backsolvable system of matrix equations}
\label{sec:solvegood}

A general formally backsolvable system of matrix equations is made
up of {\em compatibility conditions on the knowns}, equations
(\ref{tri-1}-\ref{tri-known}), and  {\em equations triangular in
the unknowns}, equations (\ref{def1:u1}-\ref{tri-last}).  One can
first verify that a completion may be possible by verifying the
compatibility conditions on the knowns
(\ref{tri-1}-\ref{tri-known}).

Next,  one attempts to use the equations triangular in the
unknowns to solve for the $u_k$ matrices. One may solve for the
(possibly non-unique) ${\bf u_{\sigma(1)}}$ using equations
(\ref{def1:u1}-\ref{defn:u1}) or determine that a solution for
${\bf u_{\sigma(1)}}$ does not exist. Notice that this may
constitute a difficult numerical problem by itself, especially if
the matrices under consideration are of large dimension. With the
results obtained for ${\bf u_{\sigma(1)}}$, one may next use
equations (\ref{def1:u2}-\ref{defn:u2}) to solve for ${\bf
u_{\sigma(2)}}$ or to determine that a solution does not exist.
This process continues until we have solved for all unknowns, that
is until
 we have
formed an inverse completion,  or have determined that a
completion is not possible. The boldface $u$'s in equations
(\ref{def1:u1}-\ref{tri-last}) indicate the unknown matrix being
solved for in each step,  if a solution can be found for each
boldface matrix.

\subsection{The strength and limitations of our method }
\label{sec:strength}

The difference in strength between Theorem~\ref{big-thm}, which was
derived automatically by computer algebra,  and Theorem~\ref{particular-thm},
which concentrated on one case and employed some human intervention,
illustrates the limitations of our automatic methods.
Theorem~\ref{particular-thm}
reduced the particular completion problem it addressed, configuration
(\ref{central-prob}),
to solving a set of compatibility conditions on the knowns and
a set of singletons defining each of the unknowns in the problem.
This is of course the
most desirable form of solution.
On the other hand, Theorem~\ref{big-thm} reduced
all but 36 of the 31,824 matrix completion problems it
addressed to essentially decoupled equations,
a highly informative but less desirable answer.

The way the stronger answer in
Theorem~\ref{particular-thm}
was derived illustrates the role of human intervention.
First, we apply Theorem~\ref{big-thm} to problem
(\ref{central-prob})
and obtain an essentially decoupled set of equations $E$.
One equation has the form
$$( a^{-1} b - d^{-1} e ) z = q$$
where $ q$ is  a polynomial.
We assumed that
$a^{-1} b - d^{-1} e $
is invertible, and implemented the assumption
by adjoining the equations defining the inverse
to $E$ and rerunning the Gr\"obner basis
 algorithm.{\footnote{
Also there was a bit of beautification of
formulas which was not essential to the form of
the result.
}}
Naturally, this solved for $z$, in
other words, produced a singleton defining $z$, and thus
derived those compatibility equations on knowns resulting from substituting
for $z$ in the equations in $E$ where it appeared.

The only human intervention
behind
Theorem~\ref{particular-thm}
was the decision that
$a^{-1} b - d^{-1} e$ is invertible.
However, a human would not know that it is critical
before Theorem~\ref{big-thm} was applied.
Without making such an invertibility assumption,
the results in Theorem~\ref{big-thm} are as far as one can go.

In fact, the invertibility assumptions and subsequent
computer manipulation, which transform
equations (\ref{quad}-\ref{lin2}) into a singleton
(\ref{z-rel2}) or (\ref{z-rel1}) in $z$, are typical of
the sort of human intervention which is required
in many problems.

\section{Solving the purely algebraic inverse matrix completion problem }

In this section,
we will describe a method for solving general matrix
completion problems of
the type described above.     The main tool we will use for
our solution of the problem is the production of a
noncommutative Gr\"{o}bner basis.
We will review Gr\"obner basis definitions and results, and
 present a pure algebra interpretation of our
matrix completion problem.   This section also contains
the formal proofs of the results presented in
Section~\ref{sec:mainresults}.

\subsection{Background on Gr\"{o}bner bases }

Gr\"{o}bner bases are a useful tool in the manipulation and analysis
of polynomial ideals.    We will review how the Gr\"{o}bner basis
may be used to

\newcounter{grobuse}
\begin{list}{\Roman{grobuse}}{\usecounter{grobuse}}
\item{Discover whether a polynomial $p$ is a member of a polynomial ideal
$I$.}
\label{grob-probs}
\item{Show two polynomial ideals $I$ and $J$  are the same,  given
generating polynomial sets $g_I$ and $g_J$.   }
\item{Transform a set of equations into an equivalent set
with a ``triangular" form, described in Section~\ref{tri-solution}.}
\end{list}

\subsubsection{Monomial orders }

Essential to the polynomial machinery we use is the existence
of a total order on the monic monomials in the polynomial
algebra under consideration.

Recall the definitions of lexicographic and graded (length)
lexicographic monomial orders on commutative monomials, as
discussed in \cite{CLS}.  The noncommutative  versions are
essentially similar,  but to ensure a well defined total order a
monomial may\footnote{In fact, this is the scheme used in the NCGB
computations.} be parsed from left to right in the tie breaking
length lexicographic order criteria.

The NCGB software uses a combination of these two types of orders,
which we will find useful. It lets one define sequential subsets
of indeterminates such that each subset is ordered with graded
lexicographic ordering within the subset,  but
in which indeterminates of a
higher set are lexicographically higher than indeterminates of a
lower set. That is, a monomial consisting of one element in a
higher set will sit higher in the monomial ordering than a
monomial consisting of the product of any number of elements in
lower sets. The NCGB notation uses the $\ll$ symbol to discern the
subset breakpoints discussed above. For example, when we write
$x_1 < x_2 < x_3 \ll x_4 $ we get that $x_1 x_2 <  x_2 x_1 $, $
x_3 x_2 x_1 < x_4 $, and $x_3 < x_1 x_2 $. We call such an
ordering {\it multigraded lexicographic}.

\subsubsection{Lead terms and Gr\"{o}bner rules }
\label{grules}

Given a polynomial $p$,  there exists a unique term of $p$ whose
monomial is highest in such an order. Denote this
$\mbox{LeadTerm}(p)$. For example, if we have $x_3 < x_2 < x_1$
then
$$ \mbox{LeadTerm}(x_1 -  x_2 x_3 + x_1^2) = x_1^2.$$

 For technical reasons, \cite{CLS},
our orders must satisfy the following relation for any polynomials $p$ and
$q$,
$$\mbox{ LeadTerm}(p) \mbox{ LeadTerm}(q) = \mbox{LeadTerm}(p \cdot q).$$
With this definition of leading term,  we can introduce
replacement rules.

 Every polynomial $p$
corresponds to a {\em replacement rule} $\Gamma(p)$,
where the
left side of the rule (LHS$\rightarrow$RHS) is the leading term
of the polynomial, and the right
side is the negative sum of the remaining
terms in the polynomial.
If we have $x_1>x_2>x_3$, then  the polynomial $x_1 -  x_2 x_3 + x_1^2$
corresponds
to the rule $x_1^2 \rightarrow x_2 x_3 - x_1$.  We may write
$$\Gamma(x_1 -  x_2 x_3 + x_1^2)  =
x_1^2 \rightarrow x_2 x_3 - x_1.$$

The next example illustrates how
we can apply a set
of replacement rules to a polynomial.

\begin{ex}
To the polynomial $$ x_1 x_1 x_2 x_1 + x_2 x_3 x_1 + x_2,$$
we may apply
$$\{x_1^2 \rightarrow x_2 x_3 - x_1, \, x_2 x_3 x_1 \rightarrow 4 x_1\}$$
to get
$$( x_2 x_3 - x_1 ) x_2 x_1 + (4 x_1) + x_2 =
x_2 x_3 x_2 x_1 -  x_1  x_2 x_1 + 4 x_1 + x_2.$$
\end{ex}

A polynomial $p$  may be {\em reduced} to a hopefully simpler form
with a set of polynomials $P$ by applying the replacement rules
$\Gamma(P)$.
Reducing a polynomial,  by a set of rules
for commutative one variable polynomials,
is similar to the
classical
Euclidean division algorithm one might use to find say
$\frac{x^3+3x^2+x+1}{x^2 + 2} $.
In the scenario of this paper, it is easier to view
reduction as a replacement scheme where rules are applied
to a polynomial to change it to a ``simpler" form.

\begin{ex}
The polynomial $x_2 x_1 x_2 x_1 + x_2$ may be reduced by the polynomial
$x_1 x_2 + 3$:

\medskip
\noindent
The Euclidean division algorithm approach:
$$ (x_2 x_1 x_2 x_1 + x_2) - x_2 (x_1 x_2 + 3) x_1  =  - 3 x_2 x_1 + x_2 $$
\medskip
\noindent
The replacement rule approach:
\begin{center}
\begin{tabular}{cccc}
 To & $ x_2 x_1 x_2 x_1 + x_2$ &  apply & $\Gamma(x_1 x_2 + 3) =
x_1 x_2 \rightarrow - 3 $
\end{tabular}
\end{center}
\begin{center}
to get $ - 3 x_2 x_1 + x_2. $
\end{center}
\end{ex}

If repeated application
of these types of rules to a polynomial transforms the equation to $0$,
then
we have shown that the polynomial under
consideration is an element of the (two-sided) ideal generated
by the relations used to create the rules.

\subsubsection{Gr\"{o}bner bases }
\label{grob-basis}

A Gr\"{o}bner basis ${\cal G}$ for a polynomial ideal $I$
enjoys the powerful
property that a polynomial $q$ is an element of the ideal if and only if
repeated application of all rules $\Gamma(g)$ arising from ${\cal G}$
send the polynomial $q$ to $0$. See \cite{CLS} for the commutative version
or \cite{HS} for the noncommutative generalization.  One says that the
Gr\"{o}bner basis solves the ideal membership problem,  which is
Problem I
given in the beginning of this section.
One may write
$$q\stackrel{ {\cal G} }{\rightarrow } 0.$$

Given generating sets  ${\cal M}$ and ${\cal N}$,  one can use this
property of Gr\"{o}bner bases to show that the ideals
$\langle {\cal M}  \rangle$ and $\langle {\cal N}  \rangle$ are
the same.
Given Gr\"{o}bner bases ${\cal G_M}$ and ${\cal G_N}$
for ${\cal M}$ and ${\cal N}$,  respectively,
we will have
\begin{equation*}
\langle {\cal M}  \rangle = \langle {\cal N}  \rangle
\end{equation*}
if
\begin{equation*}
m\stackrel{{\cal G_N}}{\rightarrow} 0 \mbox{ for all polynomials } m \in
{\cal M}
\mbox{ and }
n\stackrel{{\cal G_M}}{\rightarrow} 0 \mbox{ for all polynomials } n \in
{\cal N}.
\end{equation*}
Thus,  Problem II,  presented above can be solved, if finite
Gr\"obner bases can be found for both of the ideals.  It may be
solved, without computing a Gr\"obner basis,  by using the
reduction properties of a partial Gr\"obner basis,   as will be done in
Section~\ref{conf}.

\subsubsection{Generating a Gr\"{o}bner basis}

If the indeterminates commute,  then  a Gr\"{o}bner basis  is
always a finite set of polynomials,  and there exists an
algorithm called Buchberger's algorithm which finds this set,
given a generating set of relations for the ideal.    In the
noncommuting case, a Gr\"{o}bner basis for an ideal may be
infinite.   Nevertheless, there exists a similar algorithm due to
F. Mora \cite{FM}, which recursively defines a Gr\"{o}bner basis
and terminates if it happens to be finite.   In practice, even
this finite Gr\"{o}bner basis may be incomputable when computer
resources are taken into account.  One stops the algorithm after
a specified number of iterations, thereby generating some finite
approximation to a Gr\"{o}bner basis.

 This finite approximation,  though not exhibiting the
omniscient powers of a true Gr\"{o}bner basis,  is often useful in
reduction, as will be shown below in Sections \ref{smallb} and \ref{conf}.
These finite approximations are, of course,  made up of elements
of the ideal generated by the original relations.   There are a
few available software implementations of this noncommutative
algorithm (\cite{NCGB}, NCGB  and \cite{OPAL}, OPAL).    For our
computations,  we used the package NCGB.

The form of this (partial) Gr\"{o}bner basis is dependent
on the order under which Mora's algorithm is executed.
One places variables which one wishes eliminated
high in the order.
The order
\begin{equation*}
k_1 < k_2 < \cdots < k_l \ll u_1 \ll u_2 \ll \cdots\ll u_{2n^2-l}
\end{equation*}
will cause the output of the Gr\"{o}bner basis algorithm to
have the form of the equations (\ref{tri-1}-\ref{tri-last})
as much as is possible, as will be discussed in
Section~\ref{sec:elim}.   We may therefore be able
to solve Problem III
introduced  on page \pageref{grob-probs},  if  a good enough
approximation to a Gr\"obner basis is practically computable.

\subsubsection{Elimination theory }
\label{sec:elim}

As promised above, we now motivate the multigraded
lexicographic ordering with a
central concept in elimination theory.

\begin{defn}
A monomial order on the monic monomials in
$\mathbb{F}[x_1,\ldots,x_n]$ is said to be of the {\em $j$-th
elimination order} if a monomial containing one of
$\{x_{j+1},\ldots, x_n \}$ is ordered higher than any monomial
made up of the indeterminates $\{x_1,\ldots,x_j\}$.
\end{defn}

 A Gr\"{o}bner basis ${\cal G}$
for an ideal ${\cal I}$, created under a $j$-th elimination order,
exhibits the following ideal relation $\langle{\cal G}\bigcap
\mathbb{F}[x_1,\ldots,x_j]\rangle = {\cal I} \bigcap
\mathbb{F}[x_1,\ldots,x_j]$.    In words, this says that all
polynomials resulting from our original generating relations
which contain only the variables $\{x_1,\ldots,x_j\}$   can be
found in the ideal generated by ${\cal G}\bigcap
\mathbb{F}[x_1,\ldots,x_j]$.   This property,  well known in the
commutative case,  was extended recently to noncommutative
algebras in \cite{HS}, Theorem 11.3. We will use this property
to assist with the triangular goal described in Sections
\ref{sec:goodanswer} and \ref{sec:betteranswer}.

\subsubsection{A small basis }

Often, what a human wishes to find is not a set of relations with
the reduction properties described in Section~\ref{grob-basis},
but a small set of relations which describes the solution set. We
refer to such a small set as a small basis. Such a goal may be
accomplished by iteratively including the first $k$ elements of
our basis $B$ in a set $B_k$, creating a partial Gr\"obner basis
${\cal G}_{B_K}$ from $B_k$, and trying to reduce the other
polynomials $B\backslash {B_k}$ with this partial Gr\"{o}bner
basis  ${\cal G}_{B_k}$.   When $B_k$ has the property that the
excluded relations $B\backslash {B_k}$ are elements of the ideal
generated by the included relations $B_k$, our goal has been
achieved. {\em All} of our relations lie in the polynomial ideal
generated by $B_k$. We call such an algorithm the {\it small basis
algorithm}.

The sequence in which relations are presented to the small
basis algorithm is obviously important.   The small basis algorithm acting on
$\left(x^3,x^2,x,1\right)$ returns the unenlightening
$\left(x^3,x^2,x,1\right)$, but
when presented with $\left(1,x,x^2,x^3\right)$ the algorithm returns
$\left(1\right)$.
(Computer time required discourages some idealized implementation,
which would consider all permutations of our relations.)

\subsection{A pure algebra interpretation of
the purely algebraic partially prescribed inverse matrix
completion problem } \label{pure-alg}

This section describes
our matrix completion problem in the language of an algebraist.
The reader may skip this section,  if so desired,
with no loss of continuity to the paper.

 Labeling the known blocks as $k_i$, we may
consider the free algebra $\mathbb{F}[k_1,\ldots,k_l]$ over the
field $\mathbb{F}$ under consideration, modulo some presupposed
conditions (e.g. the invertibility of a known submatrix $k_i$,
which is expressed in ideal theoretic notation as $\langle k_i
k_i^{-1} -1, k_i^{-1} k_i - 1 \rangle$)
\begin{equation}
\label{def:S}
 S = \frac{ \mathbb{F}[k_1,\ldots,k_l]}{\langle
\mbox{ conditions on the knowns } \rangle }.
\end{equation}

Let the size of the square submatrices be $m$. Picking the known
block matrices consists of defining a map
$$\phi : S
\rightarrow M_m(\mathbb{F}).$$
Defining
\begin{equation}
\label{def:T}
 T = S[u_1,\ldots,u_{ 2 n^2 - l } ],
\end{equation}
completing the matrices $A$ and $B$ may be viewed as a map
$$\Phi : T \rightarrow M_m(\mathbb{F})$$
such that
$$\Phi \!\! \mid_{S} = \phi.$$
Here, we are interested in a special completion $\Phi$  so that
our matrices satisfy (\ref{inverse-cond}). Let {\em Flatten} be
the operation which takes a set of matrices to their constituent
blocks.  If $J$ is then the ideal generated by relations
(\ref{inverse-cond}),
\begin{equation}
\label{Jdefn}
J = \langle Flatten( A B - I , B A - I ) \rangle,
\end{equation}
our goal is achieved when $J$ lies in the kernel of $\Phi$.

Given $G$,  a generating set for $J$,  so that
$$\langle G \rangle = J,$$
our compatibility conditions on the knowns have been defined to be
$G \cap S$. These are equations (\ref{tri-1}-\ref{tri-known}) in
the backsolvable case or equations (\ref{flat-1}-\ref{endknowns})
in the decoupled case. Notice that all of the relations making up
$G \cap (J \setminus S)$ will contain at least one $u_i$.

\subsection{Nondegenerate solutions }
\label{sec:nondeg}

On closer analysis of the {\em weakly} formally backsolvable form
described in equations (\ref{tri-1}-\ref{tri-last}), the existence
of $q_{0,1}$ implies the existence of $q_{1,1},\ldots,q_{2 n^2 -l,
m_{2 n^2 -l} }$. Simply multiply $q_{0,1}$ by the appropriate
${\bf u_i}$. A more interesting set of relations has the following
property.
\begin{defn}
A set of equations $\{p_j=0: 1\le j\le s\}$ will be called {\em
\bf nondegenerate} if
\begin{equation}
\label{nonred}
p_i \not\in \langle \{ p_j \}_{j \neq i}
\rangle.
\end{equation}
\end{defn}
Condition (\ref{non-deg-cond}) in Definition~\ref{f-back}
or condition (\ref{non-deg-cond-better}) in Definition~\ref{def:ess-decoup}
is equivalent to the nondegenerate property.
Due to the infiniteness of the noncommutative Gr\"obner basis,
condition (\ref{nonred}) cannot  in general be verified. A
condition which can be verified computationally is the following.
\begin{defn}
A set of equations $\{p_j=0: 1\le j\le s\}$ will be called {\em
\bf $\ell$-nondegenerate} if
\begin{equation}
\label{k-nonred}
p_i \not\in \langle \{ p_j \}_{j \neq i}
\rangle_\ell
\end{equation}
where $\langle \{ p_j \}_{j \neq i}
\rangle_\ell$ is the $\ell$ iteration partial Gr\"obner basis created from
$\{ p_j \}_{j \neq i}$.
\end{defn}

Consider the system of equations
\begin{eqnarray}
k_1 + k_2 k_1 &=& 0 \label{knowns} \\
k_1 u_1 + k_2 k_1 u_1 &=& 0 \label{redund} \\
k_1 k_2 + k_1 u_2 k_2 u_2 + k_1 u_2 &=& 0 \label{eqn3}   \\
k_2 k_2 u_2 + u_1 k_2 k_1 &=& 0 \label{eqn4}
\end{eqnarray}
where $k_1$ and $k_2$ are known and
$u_1$ and $u_2$ are unknown.

These relations appear to be essentially decoupled,  but are only weakly
essentially decoupled,  since equation (\ref{redund}) is a
member of the ideal generated by equation  (\ref{knowns}).
By removing equation (\ref{redund}) we are left with a formally backsolvable
system of equations,  since we have an equation in one unknown $u_2$ and
an equation in two unknowns,  $u_1$ and $u_2$.  Nondegeneracy can be verified
by attempting, and failing, to reduce to $0$ each polynomial in
$\{(\ref{knowns}),(\ref{eqn3}),(\ref{eqn4})\}$
using a Gr\"obner basis created from the other two polynomials.
For this example, the three Gr\"obner bases are finite,  so this
procedure is definitive.

\subsubsection{The special case of compatibility
nondegenerate solutions }
\label{sec:nondeg-decoupled}

An approximation to nondegeneracy,  often used in this paper,
is compatibility nondegeneracy.     This form of
nondegeneracy is the condition that the equations which
contain unknowns cannot be reduced to $0$ by those equations
which contain only knowns.

\begin{defn}
We will call a weakly decoupled set of equations of the
form (\ref{flat-1}-\ref{flat-last})
{\bf compatibility nondegenerate} if
\begin{equation}
 q_{h_1,h_2} \not\in \langle \{ q_{0,1},\ldots,q_{0,m_0}
\} \rangle \mbox{ for } h_1 > 0 \mbox{ and } q_{h} \not\in \langle
\{q_{0,1},\ldots,q_{0,m_0} \} \rangle \mbox{ for } h > s_j + 2 n^2
-l -j. \label{nonred-flat}
\end{equation}
\end{defn}
That is,  the relations
(\ref{beginunk}-\ref{flat-last}) which define the $u_{\sigma(i)}$ for
$i=1,\ldots,j$ are not trivial,  and
are merely consequences of the
compatibility conditions on the knowns,
equations
(\ref{flat-1}-\ref{endknowns}).
The singleton
equations
(\ref{single}-\ref{end-singles}),  those which define
$u_{\tau(i)}$ for $i= j+1,\ldots,2 n^2 -l$,  are obviously not
trivial,  since the term $u_{\tau(i)}$ cannot be
reduced by any Gr\"obner rule containing only
$k_1,\ldots,k_l$.

We also have the computational analogue.
\begin{defn}
\label{def:k-nonred} We will call a weakly decoupled set of equations in the
form of (\ref{flat-1}-\ref{flat-last}) {\em \bf
compatibility
$\ell$-nondegenerate} if
\begin{equation}
 q_{h_1,h_2} \not\in \langle \{ q_{0,1},\ldots,q_{0,m_0}
\} \rangle_{\ell} \mbox{ for } h_1 > 0 \mbox{ and } q_{h} \not\in
\langle \{q_{0,1},\ldots,q_{0,m_0} \} \rangle_{\ell} \mbox{ for }
h
> s_j + 2 n^2 -l -j.  \label{comp-l-nonred}
\end{equation}
where  $\langle \{q_{0,1},\ldots,q_{0,m_0} \} \rangle_{\ell}$ is
the $\ell$ iteration partial Gr\"obner basis created from
$\{q_{0,1},\ldots,q_{0,m_0} \}$.
\end{defn}

A weakly formally backsolvable system of equations can also be
compatibility $\ell$-nondegenerate.
\begin{defn}
\label{def:k-nonred-fback} We will call a weakly formally backsolvable
set of equations in the
form of (\ref{tri-1}-\ref{tri-last}) {\em \bf
compatibility
$\ell$-nondegenerate} if
\begin{equation}
q_{h_1,h_2} \not\in \langle \{q_{0,1},\ldots,q_{0,m_0} \}
\rangle_{\ell} \mbox{ for } h_1 > 0   \label{comp-l-nonred-fback}
\end{equation}
where  $\langle \{q_{0,1},\ldots,q_{0,m_0} \} \rangle_{\ell}$ is
the $\ell$ iteration partial Gr\"obner basis created from
$\{q_{0,1},\ldots,q_{0,m_0}\}$.
\end{defn}

Beware that these definitions are algorithm dependent,
since the Gr\"obner basis algorithm allows for some variability
in how it is implemented.
Furthermore,  research into different
variants of the noncommutative Gr\"obner basis algorithm
has not thoroughly addressed the reduction properties of partial Gr\"obner
bases.
If one were to run the Gr\"obner basis algorithm for an infinite
number of iterations (which might result in infinitely
many polynomials),
then
one could verify
condition (\ref{nonred-flat}).

Our computational resources are, of course, finite and we can do
no such thing. Still,  our three iteration partial Gr\"obner basis
offers a computational approximation to the condition
(\ref{nonred-flat}). This form of non-redundancy given in
Definition~\ref{def:k-nonred} was used to
verify compatibility 3-nondegeneracy for all the problems in
Theorem~\ref{big-thm}, which were essentially decoupled.  All but 36
of them
were of this form.  While this is all that we did automatically on
all 31,824 cases, Theorem~\ref{particular-thm} serves to show
what one can do by further applying our Gr\"obner basis methods
to a particular case.   In that case, we gave a concise solution
to the matrix completion problem without an ``infinite
computation".

\subsection{A recipe for solving the general block matrix inverse
completion problem }
\label{sec-recipe}

We are given matrices $A$ and $B$ partitioned conformally for
matrix multiplication into $n^2$ blocks each and a configuration
of $l$ prescribed (known) and $2 n^2-l$ unknown blocks.  We may
also be given conditions on these matrices which are expressed
algebraically (e.g. invertibility,  $a a^{-1} -1 = 0$ and $a^{-1}
a -1 = 0$). We look to discover compatibility conditions on the
known matrices and formulas for the unknown matrices to solve our
problem, that is,  to ensure (\ref{inverse-cond}) is satisfied.

This paper shows
that
 this
goal may often be achieved by following the steps below.

\newcounter{probs}
\begin{list}{\Roman{probs}}{\usecounter{probs}}
\item{Fill in the known blocks of $A$ and $B$ with symbolic,
noncommuting indeterminates, $k_1,\ldots,k_l.$
}
\item{Fill in the unknown blocks of $A$ and $B$ with symbolic,
noncommuting indeterminates, $$u_1,\ldots,u_{2 n^2 - l}.$$}
\item{Create the noncommutative polynomials resulting from the
operations $A B - I$ and  $B A - I$.}
\item{Create a (noncommutative, partial) Gr\"obner basis for the polynomials
derived in step III and any assumed algebraic conditions on the
matrices under the order:
\begin{equation*}
k_1 < k_2 < \cdots < k_l \ll u_1 \ll u_2 \ll \cdots \ll u_{2 n^2 -l}
\end{equation*} }
\item{Check that the result has
some attractive form, such as those described in Section
\ref{sec:betteranswer}, (\ref{flat-1}-\ref{flat-last})
or Section \ref{sec:goodanswer}, (\ref{tri-1}-\ref{tri-last}). }
\item{ Verify that the relations defining
unknown matrices
are not merely consequences of the other relations by using
the Small Basis Algorithm or some variant of it. }
\end{list}

\noindent
The noncommutative algorithms we use are not yet well understood and,
therefore, their effectiveness on a particular class of problems
can only be determined by experimentation.

\subsection{Proof of seven unknown,  11 known theorem }
\label{big-thm-proof}

We created a {\em Mathematica} procedure,  which iteratively
searches through all permutations of seven unknown blocks and 11 known
blocks, and performs the sort of analysis described in
Section~\ref{sec-recipe}.
As described in Section \ref{configsNperms},
one may apply permutation
matrices, $\Pi$ and $\Psi$,  to $A$ and $B$  to get $\Pi^{-1} A
\Psi$ and $\Psi^{-1} B \Pi$,  and  obtain at most 36 other
equivalent configurations. This property was exploited to reduce
the computations needed from 31,824 cases to about
1,500 cases. Only one matrix inverse completion
problem was analyzed from each equivalence class.

First,  we will describe the procedures followed for a particular
configuration.   Then, we will list the pseudo-code which
performed the necessary analysis for the entire problem. Our
proof will be completed with a discussion of the results of our
{\em Mathematica} procedure.

\subsubsection{A particular configuration }

We created a two-iteration, partial Gr\"obner basis from the
polynomial matrix equations resulting from $A B$ and $B A$, along
with the invertibility relations of the knowns.

The order we used to create the Gr\"obner basis  was the
following
\begin{eqnarray}
\label{gen-order}
& & k_1 < k_2 < k_3 < k_4 < k_5
< k_6 < k_7 <
k_8 < k_9 <
k_{10} <  k_{11} \nonumber  \\
 & & \ll  u_1 \ll u_2 \ll u_3 \ll
u_4 \ll u_5  \ll u_6  \ll u_7
\end{eqnarray}
where the $k_j$ represents the $j^{th}$ known block and $u_i$
represents  the  $i^{th}$ unknown block.
Inverses have
been suppressed
in our lists of knowns for clarity.  Any listing of known variables should
be accompanied by their inverses.  These inverses are placed directly
above, and in the same group as, the original variable.   So, our order truly
begins $k_1<k_1^{-1}<k_2< k_2^{-1}< \ldots $.

The output of the Gr\"obner basis algorithm, in virtually all
cases, was of the weakly essentially decoupled form described in
Section \ref{sec:betteranswer},  equations
(\ref{flat-1}-\ref{flat-last}).

To establish weakly essentially decoupled and
compatibility 3-nondegeneracy, we used
the output of the Gr\"obner basis algorithm,  which consisted
solely of known indeterminates, equations
(\ref{flat-1}-\ref{endknowns}),  or $G\cap S$ in the language of
Section \ref{pure-alg}. We ran the Gr\"obner basis algorithm for
one more iteration on these known relations, thereby creating a three
iteration partial Gr\"obner basis. We used this partial Gr\"obner
basis to attempt to reduce the relations which contain the
unknown indeterminates.

After applying the Gr\"obner rules, associated with this partial
Gr\"obner basis, to the set of relations containing unknown
indeterminates, our set of relations still had the form given in
equations (\ref{flat-1}-\ref{flat-last}).   That is, we verified
compatibility 3-nondegeneracy as given
by condition (\ref{comp-l-nonred}).   This
verification was done by computer. This shows that the problem
associated with this particular configuration is weakly
essentially decoupled and compatibility 3-nondegenerate.

Configuration (\ref{bad}) and permutations of this configuration
were weakly formally backsolvable and compatibility 3-nondegenerate.
Theorem~\ref{big-thm} follows.

\subsubsection{Pseudo-code}

Here we give some pseudo-code with a Mathematica slant, which
performs the sort of analysis described in the above section for
all seven unknown and 11 known configurations.
An essential part of the algorithm is the function
{\tt NCMake\-GroebnerBasis[polys,k]},  which creates a $k$ iteration
partial Gr\"obner basis from $polys$.

\medskip

\noindent
{\small First, we create the relations which are implied by
the invertibility of the knowns.}

\vspace{-.5cm}

\begin{verbatim}
inverses = NCMakeRelations[{Inv, k1,k2,k3,k4,k5,k6,k7,k8,k9,k10,k11}]
\end{verbatim}

\noindent
{\small Next,  we set the monomial order for the Gr\"obner basis computation.
This order is given in (\ref{gen-order}). }

\vspace{-.5cm}

\begin{verbatim}
SetMonomialOrder[ k1<k2<k3<k4<k5<k6<k7<k8<k9<k10<k11<<u1<<u2<<u3<<u4<<u5<<u6<<u7 ]
\end{verbatim}

\noindent
{\small We then generate all permutations of seven $1$'s and eleven $0$'s.
}

\vspace{-.5cm}

\begin{verbatim}
permList = Permutations[ {0,0,0,0,0,0,0,0,0,0,0,1,1,1,1,1,1,1} ]
\end{verbatim}

\noindent
{\small
We examine all the configurations associated with the generated
permutations with the }

\vspace{-.2cm}

\noindent
{\small
following \verb+For[ ]+ loop.}

\vspace{-.5cm}

\begin{verbatim}
For[ i = 1, i++, i <= Length[ permList ],
\end{verbatim}

{\small
If a permutation (in the sense of Section~\ref{configsNperms})
of this configuration was already examined,  don't bother. }

\vspace{-.5cm}

\begin{verbatim}
   If[ MemberQ[ alreadyDoneList, permList[[i]] ]
        Continue[]
     ]
\end{verbatim}

{\small
Since no permutation of this configuration has been
 analyzed,  we add all permutations }

\vspace{-.2cm}

{\small
of this configuration to the \verb+alreadyDoneList+. }

\vspace{-.5cm}

\begin{verbatim}
   AppendTo[ alreadyDoneList, MakeTransformations[ permList[[i]] ] ]
\end{verbatim}

{\small
Next,  convert this configuration into two {\em Mathematica} matrices.}

\vspace{-.5cm}

\begin{verbatim}
   {A,B} = MakeSymbolicMatrices[ permList[[i]] ]
\end{verbatim}

{\small
Obtain the union of all relations: $A B = I$,  $B A = I$ and invertibility
of knowns. }

\vspace{-.5cm}

\begin{verbatim}
   relations = Union[ inverses, Flatten[
        MatrixMultiply[ A,B ] - IdentityMatrix[3],
        MatrixMultiply[ B,A ] - IdentityMatrix[3]  ]  ]
\end{verbatim}

{\small
Make a Gr\"obner basis from the relations generated in the
previous step.}

\vspace{-.5cm}

\begin{verbatim}
   output = NCMakeGroebnerBasis[ relations, 2 ]
\end{verbatim}

{\small Isolate the compatibility conditions on the knowns:}

\vspace{-.2cm}

{\small
 $G \cap S$ in the language
of Section~\ref{pure-alg} or equations (\ref{flat-1}-\ref{endknowns}) in
Section~\ref{sec:betteranswer}.  }

\vspace{-.5cm}

\begin{verbatim}
   polysInKnowns = FindPolysInOnlyTheVariables[ output,
                         {k1,k2,k3,k4,k5,k6,k7,k8,k9,k10,k11}]
\end{verbatim}

{\small
Reduce the output of the Gr\"obner Basis Algorithm with the compatibility
conditions found in the }

\vspace{-.2cm}

{\small previous step.  For our purposes,  this will result in a
compatibility
3-nondegenerate set of equations.}

\vspace{-.5cm}

\begin{verbatim}
   reductionSet = NCMakeGroebnerBasis[ knownPolys, 1 ]

   output = NCReduction[ output, PolyToRule[ reductionSet ]  ]
\end{verbatim}

{\small Extract unknowns, which lie in an
equation with no other unknowns, from the reduced output.
}

\vspace{-.5cm}

\begin{verbatim}
   determinedIndeterminates = PolysInOneUnknown[ output,{u1,u2,u3,u4,u5,u6,u7} ]
\end{verbatim}

{\small
Extract singleton indeterminates from the  reduced output. }

\vspace{-.5cm}

\begin{verbatim}
   singleIndeterminates = PolysExplicit[output,{u1,u2,u3,u4,u5,u6,u7}]
\end{verbatim}

{\small If all our unknown indeterminates lie in an equation in
one unknown or are singletons, then our }

\vspace{-.2cm}

{\small
solution set is of the weakly essentially decoupled
form.  Otherwise,  our solution set is not.}

\vspace{-.5cm}

\begin{verbatim}
   If[Union[determinedIndeterminates,singleIndeterminates]=={u1,u2,u3,u4,u5,u6,u7},
            Print["SUCCESSFUL"]
        ]
   Else[
            Print["UNSUCCESSFUL"]
        ]

   ]  (* End of For[] loop *)
\end{verbatim}

\subsubsection{End game }

 The problems which were strongly undetermined did not
have the formally backsolvable form.
 Of the problems which were not strongly undetermined,
there  were seven cases in which the output of the two iteration
partial Gr\"obner basis did not have the form of equations
(\ref{flat-1}-\ref{flat-last}).   For these seven cases, we
performed the same analysis,  but created a three iteration
partial Gr\"obner basis instead of halting the algorithm after
two iterations,  as was done originally.   In six of these cases,
the three iteration partial bases had the form of equations
(\ref{flat-1}-\ref{flat-last})  and were  shown to be
compatibility
3-nondegenerate.
In the case associated with configuration (\ref{bad}),
the 3-iteration partial Gr\"obner basis did not have the
essentially decoupled form, and a 4-iteration partial Gr\"obner
basis proved too difficult to compute.
Therefore,  the result stated in the theorem
follows.

The {\em Mathematica} code, associated with the pseudo-code given
above,  ran for approximately 3 days,  on a Sun Ultra II with two
166Mhz processors and 1Gb of RAM. The computer was a departmental
machine and the processes associated with these
computations were therefore given only
a portion of the total computational
resources available.   The same computations on a similar machine
dedicated to this problem might take half the time. \hfil $\Box$
\lskip

\subsection{Proof of Theorem \ref{particular-thm}}
\label{sec:part-proof}

We shall need the following lemma for our proof:

\begin{lem}[Schur]
If $x_{1,1}, x_{1,2}, x_{2,1}, x_{2,2}$ are invertible block
matrices of the same size,  then
$$
\left(
\begin{array}{cc}
x_{1,1} & x_{1,2} \\
x_{2,1} & x_{2,2}
\end{array}
\right)
$$
is invertible if and only if  $- x_{2,1}  x_{1,1}^{-1} x_{1,2} + x_{2,2}$
is invertible.
\label{schur}
\end{lem}
\begin{proof}
$$
\left(
\begin{array}{cc}
I  & 0 \\
x_{2,1} x_{1,1}^{-1} & I
\end{array}
\right)
\left(
\begin{array}{cc}
x_{1,1}  & 0 \\
0  &   - x_{2,1}  x_{1,1}^{-1} x_{1,2} + x_{2,2}
\end{array}
\right)
\left(
\begin{array}{cc}
I  & x_{1,1}^{-1} x_{1,2} \\
0 & I
\end{array}
\right)
 =
\left(
\begin{array}{cc}
x_{1,1} & x_{1,2} \\
x_{2,1} & x_{2,2}
\end{array}
\right)
$$
\end{proof}

\begin{proof}[Of Theorem]

\noindent
$\Rightarrow $\\
Creating a three iteration partial Gr\"{o}bner basis with the relations
\begin{equation}
 A B = I ,   B A = I,
\mbox{ and the invertibility of the knowns},
\label{orig-eqn}
\end{equation}
using the NCGB command
NCProcess\footnote{Appendix 2 page 7 contains
the input to the NCProcess command, the
``unraveled" equations (\ref{orig-eqn}).},
yields a
set of polynomials, which includes relations
\begin{eqnarray}
 z& =& z e z  + z d g + j b z - k c i + j a g \label{quad} \\
  a^{-1} h^{-1} - a^{-1} b j h^{-1} +a^{-1} b z i^{-1}& =&
  - d^{-1} i^{-1} - d^{-1} e  j  h^{-1}  +
d^{-1} e z i^{-1} \label{lin1}\\
 f^{-1} a^{-1} - f^{-1} g d a^{-1}  +    k^{-1} z d a^{-1}&  =&
- k^{-1} b^{-1} - f^{-1} g e b^{-1} + k^{-1} z e b^{-1} \label{lin2}
\end{eqnarray}
and relations (\ref{t-rel}-\ref{y-rel}).  See Appendix 2
pages 8-10. The order used is given on page~\pageref{order}, order
(\ref{order}). Since polynomials created through the Gr\"{o}bner
basis algorithm are in the polynomial ideal generated by the
original relations, the validity of relations
(\ref{quad}-\ref{lin2}) and (\ref{t-rel}-\ref{y-rel}) is a
consequence of relations (\ref{orig-eqn}).

For us to write the relations
(\ref{quad}-\ref{lin2}) in the form (\ref{cc-1}-\ref{z-rel2}),
we require the invertibility of $(d a^{-1} - e b^{-1})$  and $
(a^{-1} b - d^{-1} e)$.  These invertibility relations are provided
by the Schur lemma given above  since the outer matrix
(\ref{outer}),  consisting of $a$,$b$,$d$, and $e$ (all knowns),
 is assumed to be invertible.
With this, we can solve for $z$ explicitly in equations
(\ref{lin1}-\ref{lin2}) and write the relations (\ref{z-rel2}-\ref{z-rel1})
defining $z$.
Furthermore,  we may use these
definitions of $z$ to write
the relations (\ref{quad}-\ref{lin2}) as (\ref{cc-1}-\ref{cc-3}).

\noindent
$\Leftarrow $\\
The converse is again approached using a Gr\"{o}bner basis method.  As
above,  the Schur
complement formulas give the invertibility of
$(d a^{-1} - e b^{-1})$  and $(a^{-1} b - d^{-1} e)$,
  which shows that relations (\ref{quad}-\ref{lin2})
follow from (\ref{cc-1}-\ref{z-rel2}).   The question then becomes
whether or not relations (\ref{orig-eqn}) are in the ideal generated
by polynomials (\ref{t-rel}-\ref{y-rel}) and (\ref{quad}-\ref{lin2}).
We create a seven iteration partial Gr\"{o}bner basis ${\cal G}_7$  from
the polynomials (\ref{t-rel}-\ref{y-rel}) and (\ref{quad}-\ref{lin2})
with the NCGB command NCMakeGB,  under the graded (length) lexicographic
monomial order.
One can verify
that the original equations (\ref{orig-eqn}) reduce to $0$
with respect to ${\cal G}_7$.
This shows that the relations $A B = I$ and $B A = I$ are elements
of the noncommutative
polynomial ideal generated by the relations (\ref{t-rel}-\ref{y-rel}) and
(\ref{quad}-\ref{lin2}) and
the invertibility of the knowns.   The result follows.
\end{proof}

\subsection{Discovering Theorem~\ref{particular-thm} and its proof  }
\bigskip

In this section, we describe the process used to discover
our particular theorem,  Theorem~\ref{particular-thm}.
This process
follows the formal notion
of a {\em strategy}, rigorously developed in \cite{HS}.

\subsubsection{Addressing our problem }

In light of our goal,   creating polynomials in few unknowns, we
used this monomial order\footnote{ Inverses have been suppressed
in our lists of knowns for clarity.  Any listing of known
variables should be accompanied by their inverses.  These
inverses are placed directly above and in the same group as the
original variable.   So our order truly begins $a< a^{-1}<b<
b^{-1}< \ldots $.}
\begin{equation}
a<b<c<d<e<f<g<h<i<j<k\ll z \ll u \ll v \ll w \ll x \ll y   \label{order}
\end{equation}
and ran the Gr\"{o}bner basis
algorithm with an iteration limit of three.

The output of this Gr\"{o}bner basis computation included
relations (\ref{t-rel}-\ref{y-rel}),  as well as
(\ref{quad}-\ref{lin2}).
(See Appendix 2 pages 8-10 for the entire output of the GBA.)
Thus, these relations are a consequence
of the original relations.
The necessity part of the proof is complete, modulo a bit of Schur
complement beautification done in Section \ref{beaut-sect}.

\subsubsection{Converse: a smaller basis } \label{smallb}

It is true that the original relations
(\ref{orig-eqn}) are members of the ideal generated by the long and ugly
relations taking up pages 8-10 of Appendix 2 (The
{\em partiality} of a Gr\"{o}bner basis at some iteration is only in respect
to its reduction properties, and not the ideal generated by these relations.).
 We could have written these
down instead of equations (\ref{cc-1}-\ref{cc-3}), our final conclusion,
and stopped,  but
we would prefer to have a more concise set of relations which
imply the original relations.
In other words, we would like to have a smaller basis for this ideal.

The computer commands in NCGB have the ability to
simplify the basis in the manner above,  in the same step as
generating it,  by setting certain options.   However,  we did not have the
computing power,  or perhaps the patience,  to isolate the few
relations on $z$ given above (\ref{quad}-\ref{lin2}) using this
method,  under the original order. To this end,  the monomial order
was changed to graded lexicographic. In NCGB notation,  we replaced
all of the $\ll$'s with $<$'s. The graded lexicographic
order computations are often of much less computational complexity,
since monomials usually must be merely checked for number of
elements.
 When our original order
was imposed on the small basis algorithm, the two iteration application
did not complete after several days running on a Sun SPARCstation-4
computer,
while under the graded lexicographic order the algorithm
finished in a few minutes.

We tried
several different sequences, of which most gave unsatisfactory
results.  The bases found were not small enough in these cases.
  An acceptable small basis obtained through this procedure
consisted of the invertibility relations on the knowns,
the relations which give the unknowns other than $z$ in terms
of $z$ (\ref{t-rel}-\ref{y-rel}),  and relations concerning
$z$ and the knowns (\ref{quad}-\ref{lin2}).   The computer work associated
with this is given in Appendix 3.

\subsubsection{Confirmation } \label{conf}

To confirm that these relations
(\ref{quad}-\ref{lin2}, \ref{t-rel}-\ref{y-rel},
 and invertibility of the knowns) imply the original
relations,  we created a noncommutative partial Gr\"{o}bner basis from
these relations and reduced the original relations with this
partial
Gr\"{o}bner basis. The original relations all reduced to $0$.
 Thus, it was
shown that the original relations were elements of the ideal
generated by the relations given above. Interestingly enough, a
five iteration partial Gr\"{o}bner basis did not reduce the
original relations (\ref{orig-eqn}),  although a seven iteration
partial Gr\"{o}bner basis did.  (See Appendix 4.)  The order used for this
computation was again the graded lexicographic.

\subsubsection{Beautification with Schur Complements }
\label{beaut-sect}

Equations (\ref{lin1}-\ref{lin2}) are especially appreciated, because
they are linear in one unknown variable $z$.   A more satisfying situation,
though, would be to have an expression for $z$ entirely in terms of the
knowns, a singleton equation.   This may be accomplished by assuming the invertibility of
\begin{equation*}
\begin{array}{lcr}
(a^{-1} b - d^{-1} e) & \mbox{ and }&
(d a^{-1} - e b^{-1})
\end{array}
\end{equation*}
in equations (\ref{lin1}-\ref{lin2}).   By the Schur Lemma \ref{schur},
this is equivalent to the invertibility of the outer matrix
$
\left(
\begin{array}{cc}
a & b \\
d & e
\end{array}
\right),
$
since all entries of this matrix are themselves invertible.
At the outset of our investigations, we had no reason to assume this
more restrictive condition.    It was only after realizing the utility
of this assumption that we added it to our conditions.

With this, $z$ is given explicitly by the equations
(\ref{z-rel1}-\ref{z-rel2})
and each of these must satisfy the quadratic (\ref{quad}).
Hence, the equations (\ref{cc-1}-\ref{cc-3}) on the knowns
are a necessary and
sufficient set of conditions for $AB = I$ and $BA = I$.

\section{Conclusion }

In this article, we have investigated the use of noncommutative
symbolic algebra software in the analysis of  partially prescribed
inverse matrix completion problems.
We described a method for solving such problems with a computer.
We have shown that the solutions  to
all 3x3 block inverse matrix completion problems with seven unknown and
11 known blocks
are of a relatively simple form.
We presented one particular theorem, and showed how it can
be massaged into a more palatable form by making some mild assumptions
on the prescribed (known) blocks.

Finally, the author would like to express appreciation for the
effort put forth by the anonymous referee.   His or
her careful reading,  criticism,  and editing have greatly
improved this paper.

\medskip
\nocite{NCA} \nocite{MMA} \nocite{FM} \nocite{MF} \nocite{BK}
\nocite{DK}
 \nocite{GreenHK97} \nocite{DellThesis} \noindent

\bibliography{master}

\section{Appendices }

\noindent
Appendices may be found at http://arXiv.org/abs/math.LA/0101245.

\end{document}





\newtheorem{defn}{{\bf Definition}}[section]
\newtheorem{thm}{{\bf Theorem}}
\newtheorem{cor}{{\bf Corollary}}
\newtheorem{lem}{{\bf Lemma}}
\newtheorem{pro}{{\bf Proposition}}
\newtheorem{cond}{{\bf Condition}}
\newtheorem{rem}{{\bf Remark}}



\newcommand{\lskip}{\vspace{.18in}}

\newenvironment{proof}{\noindent{\bf Proof:\quad}}{\hfill $\Box$ \lskip} 
\newenvironment{ex}{\lskip \refstepcounter{defn} \noindent{\bf Example 
\thedefn}}{\lskip}
\newenvironment{fact}{\lskip \refstepcounter{defn} \noindent{\bf Fact \thedefn}}{} 






\begin{titlepage}
\title{ Appendices: Using noncommutative 
Gr\"obner bases in solving partially prescribed 
matrix inverse completion problems}

\date{ \today } 
\maketitle

\vspace{3in}
\begin{center}
{\Large 
We intend these appendices to appear on the internet for 
those interested and not be included in the publication.
}
\end{center}

\end{titlepage}

\begin{abstract}
{
We investigate  the use of noncommutative Gr\"obner bases
in solving partially prescribed matrix inverse completion
problems.     The types of problems considered here are
similar to those in [BJLW].
There the authors
gave necessary and sufficient conditions for the solution of
a two by two block matrix
completion problem.     Our approach is quite different from
theirs and relies on symbolic computer algebra.

Here we describe a general method by which all block matrix
completion problems of this type may be analyzed if sufficient computational
power is available.   We also demonstrate our method  with an analysis of
all three by three block matrix inverse completion problems
with eleven
known blocks and seven unknown.   We discover that the solutions
to all such problems are of a relatively simple form.

We then perform a more detailed analysis of a particular
problem from the
31,824  three by three block matrix
completion  problems with eleven
known blocks and seven unknown.   A solution to this problem
of the form derived in [BJLW] is presented.

Not only do we give a proof of our detailed result,
but we describe the strategy
used in discovering our theorem and proof,   since it is somewhat unusual for
these types of problems.

}
\end{abstract}


\newpage
\baselineskip=3pt
\lineskip=3pt

\voffset = -1in
\evensidemargin 0.1in
\oddsidemargin 0.1in
\textheight 9in
\textwidth 6in
\normalsize
\baselineskip=12pt
\noindent

\section{ Appendix 1 - Mathematica implementation of pseudo
code given for Theorem 1 }

\begin{verbatim}
(* First we set up needed variables and lists *)
onePerms = Permutations[{0,0,0,0,0,0,0,0,0,0,0,1,1,1,1,1,1,1}];

vars = { k1,k2,k3,k4,k5,k6,k7,k8,k9,k10,k11,u1,u2,u3,u4,u5,u6,u7 };

SetNonCommutative[ vars ];

SetNonCommutative[ Inv[k1],Inv[k2],Inv[k3],Inv[k4],Inv[k5],Inv[k6],Inv[k7],
	Inv[k8],Inv[k9],Inv[k10],Inv[k11] ];

inversesKeep =
	{-1+k1**Inv[k1],-1+k10**Inv[k10],-1+k11**Inv[k11],-1+k2**Inv[k2],
	-1+k3**Inv[k3],-1+k4**Inv[k4],-1+k5**Inv[k5],-1+k6**Inv[k6],
	-1+k7**Inv[k7],-1+k8**Inv[k8],-1+k9**Inv[k9],-1+Inv[k1]**k1,
	-1+Inv[k10]**k10,-1+Inv[k11]**k11,-1+Inv[k2]**k2,-1+Inv[k3]**k3,
	-1+Inv[k4]**k4,-1+Inv[k5]**k5,-1+Inv[k6]**k6,-1+Inv[k7]**k7,
	-1+Inv[k8]**k8,-1+Inv[k9]**k9 };

permMtcs = Permutations[ IdentityMatrix[3] ];  

invPermMtcs = Map[ Inverse, permMtcs ];

knownVars = {k1,k2,k3,k4,k5,k6,k7,k8,k9,k10,k11};
unkVars = {u1,u2,u3,u4,u5,u6,u7 };

alreadyDone = {} ;


(* Here we give some necessary functions *)

(* Make symbolic matrix from 1's notation *)
MakeSymbMatrix[ onesList_List ] := Module[{idx,unkIdx=1,knIdx=1,newMtx={} },

  For[idx=1,idx<=Length[onesList],idx++,
   
    If[ onesList[[idx]] == 1,
      AppendTo[ newMtx,unkVars[[unkIdx]] ];
      unkIdx++;,
      AppendTo[ newMtx,knownVars[[knIdx]] ];
      knIdx++;
      ];
    ];
  Return[ newMtx ];
];


(* Ask if a permutation of a matrix is in the alreadyDone list *)
PermMemberQ[  onesList_List ] :=  
	Module[ {permList={}, A,B},

  	A = Partition[onesList,9 ][[1]];
  	B = Partition[onesList,9 ][[2]];
  	A = Partition[A,3];
  	B = Partition[B,3];

  	Apermed = Map[ A.#  & ,permMtcs] ;
 
  	Aperm2 = Flatten[Table[Map[#.Apermed[[permEntry]]&, permMtcs], 
		{permEntry,1, Length[Apermed]}],1];

  	Bpermed = Map[ #.B & , invPermMtcs ];
  
  	Bperm2 = Flatten[ Table[Map[ Bpermed[[permEntry]].# &, invPermMtcs ], 
		{ permEntry,1,Length[Bpermed]}],1];
		
  	For[index =1,index<=Length[Aperm2],index++,
  		AppendTo[permList, Flatten[{Aperm2[[index]],Bperm2[[index]]}]];
  	];   
  
  	If[ Intersection[ alreadyDone,permList ] === {},
    	Return[ False ];,
    	Return[ True  ];
    ];
];
	
(* Find the unknown indeterminates which lie in equations in one unknown *)	
FindDetermined[ aList_List,currentUnknowns_List ,opts___Rule] :=
	Module[{i,len,item,vars,kn,unk,n,alllengths,rules, determList,
        allvars,allind,j,totalvars,outputToFile,fileName},

  	Clear[relations];
    determList = {};

  	rules = Union[ExpandNonCommutativeMultiply[aList]];
  	rules = PolyToRule[rules];
  	rules = Union[rules];
  	len   = Length[rules];

	Do[item = rules[[i]];

    	If[Not[item===0],
        	vars = GrabIndeterminants[item];

         	If[ Length[ Union[vars] ] == 1,
              	AppendTo[determList , vars[[1]] ];
         		];
     		];
  		,{i,1,len}];  (* End Do[] loop *)

  	Return[ determList ];
];

(* Find the indeterminates which are singletons *)
FindSingletons[ aList_List,currentUnknowns_List ,opts___Rule] :=
	Module[{i,len,item,vars,kn,unk,n,alllengths,rules, singletonList,
        allvars,allind,j,totalvars,outputToFile,fileName},

  	Clear[relations];
  	singletonList = {};

  	rules = Union[ExpandNonCommutativeMultiply[aList]];
  	rules = PolyToRule[rules];
  	rules = Union[rules];
  	len   = Length[rules];

	Do[item = rules[[i]];

    	If[Not[item===0],
       		vars = GrabIndeterminants[item];

        	(* If the head is not NCM it's a singleton !! *)
         	If[Head[item[[1]]]=!=NonCommutativeMultiply, 
				AppendTo[singletonList, item[[1]] ];
         		];
     		];
  		,{i,1,len}];  (* End Do[] loop *)

  	Return[Intersection[ singletonList, currentUnknowns ]    
  	];
];

<<Extra.TeXForm;
WriteMatrixTeX[ A_List, B_List, currNum_ ] := Module[{},

	OpenWriteForTeX["thms/mtcs"<>ToString[currNum] ];

	ExpressionToTeXFile["thms/mtcs"<>ToString[currNum]<>".tex" ,
		OutputAMatrix[A] ];

	ExpressionToTeXFile["thms/mtcs"<>ToString[currNum]<>".tex" ,
		OutputAMatrix[B] ];

	WriteString["thms/mtcs"<>ToString[currNum]<>".tex",
		"\\end{document}" ];

	Close["thms/mtcs"<>ToString[currNum]<>".tex"];
]; 

OpenWrite["AnswersForCompleteMtx"]; 

undetermList = 
 {683,684,689,695,706,708,710,719,725,729,748,749,752,755,760,769,779,784,1656,
  1675,1701,2580,2599,2625,4296,4315,4341,5957,5958,5963,5969,5980,5982,5984,
  5993,5999,6003,6022,6023,6026,6029,6034,6043,6053,6058,6090,6109,6135,6174,
  6193,6219,6294,6313,6339,7299,7318,7344,7486,7487,7491,7498,7499,7500,7501,
  7520,7521,7540,7556,7575,7606,7607,7608,7609,7610,7611,7738,7757,7821,7948,
  7967,8031,9070,9089,9153,10357,10376,10440,11244,11263,11289,12304,12323,
  12349,13409,13410,13414,13421,13422,13423,13424,13443,13444,13463,13479,
  13498,13529,13530,13531,13532,13533,13534,13535,13554,13618,13829,13848,
  13912,14159,14178,14242,15941,15960,16024,16249,16268,16294,17364,17383,
  17447,17943,17962,18026,20312,20331,20357,22914,22915,22919,22926,22927,
  22928,22929,22948,22949,22968,22984,23003,23034,23035,23036,23037,23038,
  23039,23040,23059,23123,23124,23143,23207,23574,23593,23657,24069,24088,
  24152,24257,24276,24302,25372,25391,25455,26786,26805,26869,28954,28973,
  29037,29789,29808,29872,30476,30477,30482,30488,30499,30501,30503,30512,
  30518,30522,30541,30542,30545,30548,30553,30562,30572,30577,30609,30628,
  30654,30693,30712,30738,30813,30832,30858,30978,30997,31023,31198,31217,
  31243,31484,31503,31529} ;


(* Here is the MAIN LOOP *)
For[ currIndex = 1 , currIndex <= Length[onePerms],
	currIndex ++,
	
	Write["AnswersForCompleteMtx", "permutation#",
	      currIndex ];   
	
	If[ PermMemberQ[ onePerms[[currIndex]] ],
	    Print[ "Found One Already Done !!!! " ];
	    Write["AnswersForCompleteMtx",
	      " Found one already done", 
	      onePerms[[currIndex]] ];
	    
		Continue[];
		];

	(*  Only add to alreadyDone list if a perm is not in it *)
	AppendTo[ alreadyDone, 
      	onePerms[[currIndex]] ];

	If[ MemberQ[undetermList, currIndex  ],
		Write["AnswersForCompleteMtx",
          "This problem is of UNDETERMINED FORM. ",
          onePerms[[currIndex]] ];
		Continue[];
		];

    alreadyDone = Union[ alreadyDone ];
	
	newMtcs = MakeSymbMatrix[onePerms[[currIndex]] ]; 
	
	mtrcs = Partition[ newMtcs ,9 ];	
	
	matrixA = Partition[ mtrcs[[1]] ,3 ];
	matrixB = Partition[ mtrcs[[2]] ,3 ];

	WriteMatrixTeX[matrixA, matrixB, currIndex ];

	oneway = MatMult[matrixA,matrixB] - IdentityMatrix[3];

	otherway = MatMult[matrixB,matrixA] - IdentityMatrix[3]; 
	
	start = Flatten[{ oneway, otherway  }];

	inverses = inversesKeep;

	start=Join[start,inverses];

	ClearMonomialOrderAll[];
 	SetMonomialOrder[ {{k1,Inv[k1],k2,Inv[k2],k3,Inv[k3],k4,Inv[k4],
		k5,Inv[k5],k6,Inv[k6],k7,Inv[k7],k8,Inv[k8],k9,Inv[k9],
		k10,Inv[k10],k11,Inv[k11]},{u1},{u2},{u3},{u4},{u5},{u6},{u7}}];

	fileName = "thms/WoerdOutput-"<>ToString[currIndex];

    (* This function creates a Noncommutative Groebner basis and
		a TeX file describing the output *)
	answer = NCProcess[start,2,2,1,1,
     	fileName ,RR->True, SB->True, SBByCat->False, NCCV->False 
		] ; 

	If[  Complement[unkVars, Union[
 		FindDetermined[answer[[3]], unkVars],
		FindSingletons[ answer[[3]], unkVars] ] ] =!= {} ,
 
		Write["AnswersForCompleteMtx", "Didn't WORK correctly.",
	   		onePerms[[currIndex]]    ];,
		(* else *)

		Write["AnswersForCompleteMtx", "Worked fine.",
	   		onePerms[[currIndex]]    ];  
		
    	];   (*  end If[]  *)
    
	Open["thms/outRels"<>ToString[currIndex]<>".m"];
	Write["thms/outRels"<>ToString[currIndex]<>".m",
	 	answer[[3]] ];
	Close["thms/outRels"<>ToString[currIndex]<>".m"];	  

]; (* end for *) 

Close[ "AnswersForCompleteMtx"];







(**********  Redefine SmallBasis[ ] for decoupled analysis *****)

SmallBasis[input_List,keepListInverses_List,iterationCount_Integer]:=

Module[ {},

     singlePolys =  FindSinglePolys[ result, { u1,u2,u3,u4,
          u5,u6,u7 } ] ;

     result = Complement [ result,   singlePolys ];

     keepList = inversesKeep ;

     ClearMonomialOrderAll[];
     SetMonomialOrder[ { k1,Inv[k1],k2,Inv[k2],k3,Inv[k3],k4,
          Inv[k4],k5,Inv[k5],k6,Inv[k6],k7,Inv[k7], 
          k8,Inv[k8],k9,Inv[k9],k10,Inv[k10],k11,Inv[k11],
          u1,u2,u3,u4,u5,u6,u7}];     
         
     knownPolys = FindKnownPolys[ result, { u1,u2,u3,u4,
          u5,u6,u7 } ] ;
         
     (* Here we order the polys since order matters for the Small Basis Alg *)
     result = Flatten[ { knownPolys,
          FindDeterminedPolys[ result, { u1,u2,u3,u4,
               u5,u6,u7 } ], 
               Complement [ result, Union[  knownPolys ,
                    FindDeterminedPolys[ result, { u1,u2,u3,u4,
                          u5,u6,u7 } ] 
                          ] 
                ] } 
           ];
         
     keepList = Union[ keepList, knownPolys ];
         
     wRules = PolyToRule[ NCMakeGB[ keepList, 1 ] ];
         
     result = Reduction[ Complement[ result , keepList ], wRules ];
         
     result = Union[ result, keepList, singlePolys  ];
                 
     (* Restore normalcy for Regular Output to look nice *)
     result = Union[ result, singlePolys, keepList  ];
         
     ClearMonomialOrderAll[];
     SetMonomialOrder[ {k1,Inv[k1],k2,Inv[k2],k3,Inv[k3],k4,
          Inv[k4], k5,Inv[k5],
          k6,Inv[k6],k7,Inv[k7],k8,Inv[k8],k9,Inv[k9],k10,Inv[k10],
          k11,Inv[k11]},{u1},{u2},{u3},{u4},{u5},{u6},{u7} ];
     
];





\end{verbatim}

\section { Appendix 2  - The First Run for Theorem 2 }
\label{sec-attach1}

\noindent
Input = 
$
\\-1 + a\,
 x + b\,
 j + t\,
 h\\
a\,
 f + b\,
 k + t\,
 y\\
a\,
 g + b\,
 z + t\,
 i\\
c\,
 h + u\,
 x + v\,
 j\\
-1 + c\,
 y + u\,
 f + v\,
 k\\
c\,
 i + u\,
 g + v\,
 z\\
d\,
 x + e\,
 j + w\,
 h\\
d\,
 f + e\,
 k + w\,
 y\\
-1 + d\,
 g + e\,
 z + w\,
 i\\
-1 + f\,
 u + g\,
 d + x\,
 a\\
f\,
 c + g\,
 w + x\,
 t\\
f\,
 v + g\,
 e + x\,
 b\\
h\,
 a + i\,
 d + y\,
 u\\
-1 + h\,
 t + i\,
 w + y\,
 c\\
h\,
 b + i\,
 e + y\,
 v\\
j\,
 a + k\,
 u + z\,
 d\\
j\,
 t + k\,
 c + z\,
 w\\
-1 + j\,
 b + k\,
 v + z\,
 e\\
a\,
 a^{ -1 } == 1\\
a^{ -1 }\,
 a == 1\\
b\,
 b^{ -1 } == 1\\
b^{ -1 }\,
 b == 1\\
c\,
 c^{ -1 } == 1\\
c^{ -1 }\,
 c == 1\\
d\,
 d^{ -1 } == 1\\
d^{ -1 }\,
 d == 1\\
e\,
 e^{ -1 } == 1\\
e^{ -1 }\,
 e == 1\\
f\,
 f^{ -1 } == 1\\
f^{ -1 }\,
 f == 1\\
g\,
 g^{ -1 } == 1\\
g^{ -1 }\,
 g == 1\\
h\,
 h^{ -1 } == 1\\
h^{ -1 }\,
 h == 1\\
i\,
 i^{ -1 } == 1\\
i^{ -1 }\,
 i == 1\\
j\,
 j^{ -1 } == 1\\
j^{ -1 }\,
 j == 1\\
k\,
 k^{ -1 } == 1\\
k^{ -1 }\,
 k == 1\\
$
File Name = MatrixInverseAnswer-3\\
NCMakeGB Iterations = 2\\
NCMakeGB on Digested Iterations = 3\\
SmallBasis Iterations = 3\\
SmallBasis on Knowns Iterations = 4\\
Deselect$\rightarrow \{\}$\\
UserSelect$\rightarrow \{\}$\\
RR$\rightarrow $True\\
RRByCat$\rightarrow $False\\
SB$\rightarrow $False\\
SBByCat$\rightarrow $False\\
DegreeCap$\rightarrow $12\\
DegreeSumCap$\rightarrow $80\\
DegreeCapSB$\rightarrow $13\\
DegreeSumCapSB$\rightarrow $81\\
NCCV$\rightarrow $False\\
THE ORDER IS NOW THE FOLLOWING:\hfil\break
$
a<a^{ -1 }<b<b^{ -1 }<c<c^{ -1 }<d<d^{ -1 }<e<e^{ -1 }<f<f^{ -1 }<g<g^{ -1 }<h<h^{ -1 }<i<i^{ -1 }<j<j^{ -1 }<k<k^{ -1 }\ll
z\ll
y<x<t<u<v<w\\
$
\label{start-digested}
\rule[2pt]{6in}{4pt}\hfil\break
\rule[2pt]{1.879in}{4pt}
\ YOUR SESSION HAS DIGESTED\ 
\rule[2pt]{1.879in}{4pt}\hfil\break
\rule[2pt]{1.923in}{4pt}
\ THE FOLLOWING RELATIONS\ 
\rule[2pt]{1.923in}{4pt}\hfil\break
\rule[2pt]{6in}{4pt}\hfil\break
THE FOLLOWING VARIABLES HAVE BEEN SOLVED FOR:\hfil\break
$\{t,
$ $
u,
$ $
v,
$ $
w,
$ $
x,
$ $
y\}$
\smallskip\\
The corresponding rules are the following:\smallskip\\
\begin{minipage}{6in}
$
t\rightarrow -1\,
 a\,
 g\,
 i^{ -1 } - b\,
 z\,
 i^{ -1 }
$
\end{minipage}\medskip\\
\begin{minipage}{6in}
$
u\rightarrow -1\,
 k^{ -1 }\,
 j\,
 a - k^{ -1 }\,
 z\,
 d
$
\end{minipage}\medskip\\
\begin{minipage}{6in}
$
v\rightarrow k^{ -1 } - k^{ -1 }\,
 j\,
 b - k^{ -1 }\,
 z\,
 e
$
\end{minipage}\medskip\\
\begin{minipage}{6in}
$
w\rightarrow i^{ -1 } - d\,
 g\,
 i^{ -1 } - e\,
 z\,
 i^{ -1 }
$
\end{minipage}\medskip\\
\begin{minipage}{6in}
$
x\rightarrow a^{ -1 } + f\,
 k^{ -1 }\,
 j - g\,
 d\,
 a^{ -1 } + f\,
 k^{ -1 }\,
 z\,
 d\,
 a^{ -1 }
$
\end{minipage}\medskip\\
\begin{minipage}{6in}
$
y\rightarrow c^{ -1 }\,
 k^{ -1 }\,
 j\,
 a\,
 f + c^{ -1 }\,
 k^{ -1 }\,
 j\,
 b\,
 k + c^{ -1 }\,
 k^{ -1 }\,
 z\,
 d\,
 f + c^{ -1 }\,
 k^{ -1 }\,
 z\,
 e\,
 k
$
\end{minipage}\medskip\\
\rule[3pt]{6in}{.7pt}\\
The expressions with unknown variables $\{\}$\\
and knowns $\{a,
$ $
b,
$ $
c,
$ $
d,
$ $
e,
$ $
f,
$ $
g,
$ $
h,
$ $
i,
$ $
j,
$ $
k,
$ $
a^{ -1 },
$ $
b^{ -1 },
$ $
c^{ -1 },
$ $
d^{ -1 },
$ $
e^{ -1 },
$ $
f^{ -1 },
$ $
g^{ -1 },
$ $
h^{ -1 },
$ $
i^{ -1 },
$ $
j^{ -1 },
$ $
k^{ -1 }\}$\smallskip\\
\begin{minipage}{6in}
$
a\,
 a^{ -1 }\rightarrow 1
$
\end{minipage}\medskip \\
\begin{minipage}{6in}
$
b\,
 b^{ -1 }\rightarrow 1
$
\end{minipage}\medskip \\
\begin{minipage}{6in}
$
c\,
 c^{ -1 }\rightarrow 1
$
\end{minipage}\medskip \\
\begin{minipage}{6in}
$
d\,
 d^{ -1 }\rightarrow 1
$
\end{minipage}\medskip \\
\begin{minipage}{6in}
$
e\,
 e^{ -1 }\rightarrow 1
$
\end{minipage}\medskip \\
\begin{minipage}{6in}
$
f\,
 f^{ -1 }\rightarrow 1
$
\end{minipage}\medskip \\
\begin{minipage}{6in}
$
g\,
 g^{ -1 }\rightarrow 1
$
\end{minipage}\medskip \\
\begin{minipage}{6in}
$
h\,
 h^{ -1 }\rightarrow 1
$
\end{minipage}\medskip \\
\begin{minipage}{6in}
$
i\,
 i^{ -1 }\rightarrow 1
$
\end{minipage}\medskip \\
\begin{minipage}{6in}
$
j\,
 j^{ -1 }\rightarrow 1
$
\end{minipage}\medskip \\
\begin{minipage}{6in}
$
k\,
 k^{ -1 }\rightarrow 1
$
\end{minipage}\medskip \\
\begin{minipage}{6in}
$
a^{ -1 }\,
 a\rightarrow 1
$
\end{minipage}\medskip \\
\begin{minipage}{6in}
$
b^{ -1 }\,
 b\rightarrow 1
$
\end{minipage}\medskip \\
\begin{minipage}{6in}
$
c^{ -1 }\,
 c\rightarrow 1
$
\end{minipage}\medskip \\
\begin{minipage}{6in}
$
d^{ -1 }\,
 d\rightarrow 1
$
\end{minipage}\medskip \\
\begin{minipage}{6in}
$
e^{ -1 }\,
 e\rightarrow 1
$
\end{minipage}\medskip \\
\begin{minipage}{6in}
$
f^{ -1 }\,
 f\rightarrow 1
$
\end{minipage}\medskip \\
\begin{minipage}{6in}
$
g^{ -1 }\,
 g\rightarrow 1
$
\end{minipage}\medskip \\
\begin{minipage}{6in}
$
h^{ -1 }\,
 h\rightarrow 1
$
\end{minipage}\medskip \\
\begin{minipage}{6in}
$
i^{ -1 }\,
 i\rightarrow 1
$
\end{minipage}\medskip \\
\begin{minipage}{6in}
$
j^{ -1 }\,
 j\rightarrow 1
$
\end{minipage}\medskip \\
\begin{minipage}{6in}
$
k^{ -1 }\,
 k\rightarrow 1
$
\end{minipage}\\
\rule[2pt]{6in}{1pt}\hfil\break
\rule[2.5pt]{1.701in}{1pt}
\ USER CREATIONS APPEAR BELOW\ 
\rule[2.5pt]{1.701in}{1pt}\hfil\break
\rule[2pt]{6in}{1pt}\hfil\break
\rule[2pt]{6in}{4pt}\hfil\break
\rule[2pt]{1.45in}{4pt}
\ SOME RELATIONS WHICH APPEAR BELOW\ 
\rule[2pt]{1.45in}{4pt}\hfil\break
\rule[2pt]{2.18in}{4pt}
\ MAY BE UNDIGESTED\ 
\rule[2pt]{2.18in}{4pt}\hfil\break
\rule[2pt]{6in}{4pt}\hfil\break
THE FOLLOWING VARIABLES HAVE NOT BEEN SOLVED FOR:\hfil\break
$\{a,
$ $
b,
$ $
c,
$ $
d,
$ $
e,
$ $
f,
$ $
g,
$ $
h,
$ $
i,
$ $
j,
$ $
k,
$ $
z,
$ $
a^{ -1 },
$ $
b^{ -1 },
$ $
c^{ -1 },
$ $
d^{ -1 },
$ $
e^{ -1 },
$ $
f^{ -1 },
$ $
g^{ -1 },
$ $
h^{ -1 },
$ $
i^{ -1 },
$ $
j^{ -1 },
$ $
k^{ -1 }\}$
\smallskip\\
\rule[3pt]{6in}{.7pt}\\
The expressions with unknown variables $\{z\}$\\
and knowns $\{a,
$ $
b,
$ $
c,
$ $
d,
$ $
e,
$ $
f,
$ $
g,
$ $
h,
$ $
i,
$ $
j,
$ $
k,
$ $
a^{ -1 },
$ $
b^{ -1 },
$ $
c^{ -1 },
$ $
d^{ -1 },
$ $
e^{ -1 },
$ $
f^{ -1 },
$ $
h^{ -1 },
$ $
i^{ -1 },
$ $
k^{ -1 }\}$\smallskip\\
\begin{minipage}{6in}
$
z\,
 e\,
 z\rightarrow z - j\,
 a\,
 g - j\,
 b\,
 z + k\,
 c\,
 i - z\,
 d\,
 g
$
\end{minipage}\medskip \\
\begin{minipage}{6in}
$
d^{ -1 }\,
 e\,
 z\,
 i^{ -1 }\rightarrow a^{ -1 }\,
 h^{ -1 } + d^{ -1 }\,
 i^{ -1 } - a^{ -1 }\,
 b\,
 j\,
 h^{ -1 } + a^{ -1 }\,
 b\,
 z\,
 i^{ -1 } + d^{ -1 }\,
 e\,
 j\,
 h^{ -1 }
$
\end{minipage}\medskip \\
\begin{minipage}{6in}
$
k^{ -1 }\,
 z\,
 e\,
 b^{ -1 }\rightarrow f^{ -1 }\,
 a^{ -1 } + k^{ -1 }\,
 b^{ -1 } - f^{ -1 }\,
 g\,
 d\,
 a^{ -1 } + f^{ -1 }\,
 g\,
 e\,
 b^{ -1 } + k^{ -1 }\,
 z\,
 d\,
 a^{ -1 }
$
\end{minipage}\medskip \\
\begin{minipage}{6in}
$
i\,
 e\,
 z\,
 i^{ -1 }\,
 c^{ -1 }\rightarrow -1\,
 h\,
 a\,
 g\,
 i^{ -1 }\,
 c^{ -1 } - h\,
 b\,
 z\,
 i^{ -1 }\,
 c^{ -1 } - i\,
 d\,
 g\,
 i^{ -1 }\,
 c^{ -1 } + c^{ -1 }\,
 k^{ -1 }\,
 j\,
 a\,
 f + c^{ -1 }\,
 k^{ -1 }\,
 j\,
 b\,
 k + c^{ -1 }\,
 k^{ -1 }\,
 z\,
 d\,
 f + c^{ -1 }\,
 k^{ -1 }\,
 z\,
 e\,
 k
$
\end{minipage}\medskip \\
\begin{minipage}{6in}
$
k^{ -1 }\,
 z\,
 d\,
 a^{ -1 }\,
 b\rightarrow -1\,
 k^{ -1 } - f^{ -1 }\,
 g\,
 e - f^{ -1 }\,
 a^{ -1 }\,
 b + k^{ -1 }\,
 z\,
 e + f^{ -1 }\,
 g\,
 d\,
 a^{ -1 }\,
 b
$
\end{minipage}\medskip \\
\begin{minipage}{6in}
$
k^{ -1 }\,
 z\,
 e\,
 k\,
 c\rightarrow c\,
 h\,
 a\,
 g\,
 i^{ -1 } + c\,
 h\,
 b\,
 z\,
 i^{ -1 } + c\,
 i\,
 d\,
 g\,
 i^{ -1 } + c\,
 i\,
 e\,
 z\,
 i^{ -1 } - k^{ -1 }\,
 j\,
 a\,
 f\,
 c - k^{ -1 }\,
 j\,
 b\,
 k\,
 c - k^{ -1 }\,
 z\,
 d\,
 f\,
 c
$
\end{minipage}\medskip \\
\begin{minipage}{6in}
$
a\,
 f\,
 k^{ -1 }\,
 z\,
 d\,
 a^{ -1 }\rightarrow -1\,
 b\,
 j - a\,
 f\,
 k^{ -1 }\,
 j + a\,
 g\,
 d\,
 a^{ -1 } + a\,
 g\,
 i^{ -1 }\,
 h + b\,
 z\,
 i^{ -1 }\,
 h
$
\end{minipage}\medskip \\
\begin{minipage}{6in}
$
d\,
 f\,
 k^{ -1 }\,
 z\,
 d\,
 a^{ -1 }\rightarrow -1\,
 d\,
 a^{ -1 } - e\,
 j - i^{ -1 }\,
 h - d\,
 f\,
 k^{ -1 }\,
 j + d\,
 g\,
 d\,
 a^{ -1 } + d\,
 g\,
 i^{ -1 }\,
 h + e\,
 z\,
 i^{ -1 }\,
 h
$
\end{minipage}\medskip \\
\begin{minipage}{6in}
$
f\,
 k^{ -1 }\,
 z\,
 d\,
 a^{ -1 }\,
 h^{ -1 }\rightarrow g\,
 i^{ -1 } - f\,
 k^{ -1 }\,
 j\,
 h^{ -1 } + g\,
 d\,
 a^{ -1 }\,
 h^{ -1 } - a^{ -1 }\,
 b\,
 j\,
 h^{ -1 } + a^{ -1 }\,
 b\,
 z\,
 i^{ -1 }
$
\end{minipage}\medskip \\
\begin{minipage}{6in}
$
i\,
 e\,
 z\,
 i^{ -1 }\,
 h\,
 a\rightarrow h\,
 a + i\,
 d - h\,
 a\,
 g\,
 d + h\,
 b\,
 j\,
 a - i\,
 d\,
 g\,
 d + i\,
 e\,
 j\,
 a + h\,
 a\,
 f\,
 k^{ -1 }\,
 j\,
 a + h\,
 a\,
 f\,
 k^{ -1 }\,
 z\,
 d - h\,
 a\,
 g\,
 i^{ -1 }\,
 h\,
 a - h\,
 b\,
 z\,
 i^{ -1 }\,
 h\,
 a + i\,
 d\,
 f\,
 k^{ -1 }\,
 j\,
 a + i\,
 d\,
 f\,
 k^{ -1 }\,
 z\,
 d - i\,
 d\,
 g\,
 i^{ -1 }\,
 h\,
 a
$
\end{minipage}\medskip \\
\begin{minipage}{6in}
$
i\,
 e\,
 z\,
 i^{ -1 }\,
 h\,
 b\rightarrow -1\,
 h\,
 a\,
 f\,
 k^{ -1 } - h\,
 a\,
 g\,
 e + h\,
 b\,
 j\,
 b - i\,
 d\,
 f\,
 k^{ -1 } - i\,
 d\,
 g\,
 e + i\,
 e\,
 j\,
 b + h\,
 a\,
 f\,
 k^{ -1 }\,
 j\,
 b + h\,
 a\,
 f\,
 k^{ -1 }\,
 z\,
 e - h\,
 a\,
 g\,
 i^{ -1 }\,
 h\,
 b - h\,
 b\,
 z\,
 i^{ -1 }\,
 h\,
 b + i\,
 d\,
 f\,
 k^{ -1 }\,
 j\,
 b + i\,
 d\,
 f\,
 k^{ -1 }\,
 z\,
 e - i\,
 d\,
 g\,
 i^{ -1 }\,
 h\,
 b
$
\end{minipage}\medskip \\
\begin{minipage}{6in}
$
k\,
 c\,
 i\,
 e\,
 z\,
 i^{ -1 }\rightarrow j\,
 a\,
 f\,
 c + j\,
 b\,
 k\,
 c + z\,
 d\,
 f\,
 c + z\,
 e\,
 k\,
 c - k\,
 c\,
 h\,
 a\,
 g\,
 i^{ -1 } - k\,
 c\,
 h\,
 b\,
 z\,
 i^{ -1 } - k\,
 c\,
 i\,
 d\,
 g\,
 i^{ -1 }
$
\end{minipage}\medskip \\
\begin{minipage}{6in}
$
z\,
 d\,
 f\,
 k^{ -1 }\,
 z\,
 d\rightarrow -1\,
 z\,
 d + j\,
 a\,
 g\,
 d - j\,
 b\,
 j\,
 a + k\,
 c\,
 h\,
 a + z\,
 d\,
 g\,
 d - z\,
 e\,
 j\,
 a - j\,
 a\,
 f\,
 k^{ -1 }\,
 j\,
 a - j\,
 a\,
 f\,
 k^{ -1 }\,
 z\,
 d - z\,
 d\,
 f\,
 k^{ -1 }\,
 j\,
 a
$
\end{minipage}\medskip \\
\begin{minipage}{6in}
$
z\,
 d\,
 f\,
 k^{ -1 }\,
 z\,
 e\rightarrow j\,
 b + j\,
 a\,
 f\,
 k^{ -1 } + j\,
 a\,
 g\,
 e - j\,
 b\,
 j\,
 b + k\,
 c\,
 h\,
 b + z\,
 d\,
 f\,
 k^{ -1 } + z\,
 d\,
 g\,
 e - z\,
 e\,
 j\,
 b - j\,
 a\,
 f\,
 k^{ -1 }\,
 j\,
 b - j\,
 a\,
 f\,
 k^{ -1 }\,
 z\,
 e - z\,
 d\,
 f\,
 k^{ -1 }\,
 j\,
 b
$
\end{minipage}\medskip \\
\begin{minipage}{6in}
$
k^{ -1 }\,
 z\,
 d\,
 f\,
 k^{ -1 }\,
 z\rightarrow -1\,
 k^{ -1 }\,
 z + c\,
 h\,
 a\,
 d^{ -1 } + k^{ -1 }\,
 j\,
 a\,
 g + k^{ -1 }\,
 z\,
 d\,
 g - k^{ -1 }\,
 j\,
 a\,
 f\,
 k^{ -1 }\,
 z - k^{ -1 }\,
 j\,
 b\,
 j\,
 a\,
 d^{ -1 } - k^{ -1 }\,
 z\,
 e\,
 j\,
 a\,
 d^{ -1 } - k^{ -1 }\,
 j\,
 a\,
 f\,
 k^{ -1 }\,
 j\,
 a\,
 d^{ -1 } - k^{ -1 }\,
 z\,
 d\,
 f\,
 k^{ -1 }\,
 j\,
 a\,
 d^{ -1 }
$
\end{minipage}\medskip \\
\begin{minipage}{6in}
$
b\,
 z\,
 i^{ -1 }\,
 h\,
 b\,
 z\,
 i^{ -1 }\rightarrow a\,
 f\,
 c - b\,
 z\,
 i^{ -1 } - a\,
 g\,
 d\,
 g\,
 i^{ -1 } - a\,
 g\,
 e\,
 z\,
 i^{ -1 } + b\,
 j\,
 a\,
 g\,
 i^{ -1 } + b\,
 j\,
 b\,
 z\,
 i^{ -1 } - a\,
 g\,
 i^{ -1 }\,
 h\,
 a\,
 g\,
 i^{ -1 } - a\,
 g\,
 i^{ -1 }\,
 h\,
 b\,
 z\,
 i^{ -1 } - b\,
 z\,
 i^{ -1 }\,
 h\,
 a\,
 g\,
 i^{ -1 }
$
\end{minipage}\medskip \\
\begin{minipage}{6in}
$
e\,
 z\,
 i^{ -1 }\,
 h\,
 b\,
 z\,
 i^{ -1 }\rightarrow d\,
 f\,
 c + d\,
 g\,
 i^{ -1 } - d\,
 g\,
 d\,
 g\,
 i^{ -1 } - d\,
 g\,
 e\,
 z\,
 i^{ -1 } + e\,
 j\,
 a\,
 g\,
 i^{ -1 } + e\,
 j\,
 b\,
 z\,
 i^{ -1 } + i^{ -1 }\,
 h\,
 a\,
 g\,
 i^{ -1 } + i^{ -1 }\,
 h\,
 b\,
 z\,
 i^{ -1 } - d\,
 g\,
 i^{ -1 }\,
 h\,
 a\,
 g\,
 i^{ -1 } - d\,
 g\,
 i^{ -1 }\,
 h\,
 b\,
 z\,
 i^{ -1 } - e\,
 z\,
 i^{ -1 }\,
 h\,
 a\,
 g\,
 i^{ -1 }
$
\end{minipage}\medskip \\
\begin{minipage}{6in}
$
z\,
 i^{ -1 }\,
 h\,
 b\,
 z\,
 i^{ -1 }\,
 c^{ -1 }\rightarrow -1\,
 z\,
 i^{ -1 }\,
 c^{ -1 } + b^{ -1 }\,
 a\,
 f + j\,
 a\,
 g\,
 i^{ -1 }\,
 c^{ -1 } + j\,
 b\,
 z\,
 i^{ -1 }\,
 c^{ -1 } - z\,
 i^{ -1 }\,
 h\,
 a\,
 g\,
 i^{ -1 }\,
 c^{ -1 } - b^{ -1 }\,
 a\,
 g\,
 d\,
 g\,
 i^{ -1 }\,
 c^{ -1 } - b^{ -1 }\,
 a\,
 g\,
 e\,
 z\,
 i^{ -1 }\,
 c^{ -1 } - b^{ -1 }\,
 a\,
 g\,
 i^{ -1 }\,
 h\,
 a\,
 g\,
 i^{ -1 }\,
 c^{ -1 } - b^{ -1 }\,
 a\,
 g\,
 i^{ -1 }\,
 h\,
 b\,
 z\,
 i^{ -1 }\,
 c^{ -1 }
$
\end{minipage}\medskip \\
\begin{minipage}{6in}
$
k^{ -1 }\,
 z\,
 d\,
 f\,
 k^{ -1 }\,
 j\,
 b\,
 e^{ -1 }\rightarrow k^{ -1 }\,
 z - c\,
 h\,
 a\,
 d^{ -1 } + c\,
 h\,
 b\,
 e^{ -1 } + k^{ -1 }\,
 j\,
 b\,
 e^{ -1 } + k^{ -1 }\,
 j\,
 a\,
 f\,
 k^{ -1 }\,
 e^{ -1 } + k^{ -1 }\,
 j\,
 b\,
 j\,
 a\,
 d^{ -1 } - k^{ -1 }\,
 j\,
 b\,
 j\,
 b\,
 e^{ -1 } + k^{ -1 }\,
 z\,
 d\,
 f\,
 k^{ -1 }\,
 e^{ -1 } + k^{ -1 }\,
 z\,
 e\,
 j\,
 a\,
 d^{ -1 } - k^{ -1 }\,
 z\,
 e\,
 j\,
 b\,
 e^{ -1 } + k^{ -1 }\,
 j\,
 a\,
 f\,
 k^{ -1 }\,
 j\,
 a\,
 d^{ -1 } - k^{ -1 }\,
 j\,
 a\,
 f\,
 k^{ -1 }\,
 j\,
 b\,
 e^{ -1 } + k^{ -1 }\,
 z\,
 d\,
 f\,
 k^{ -1 }\,
 j\,
 a\,
 d^{ -1 }
$
\end{minipage}\medskip \\
\begin{minipage}{6in}
$
e^{ -1 }\,
 d\,
 g\,
 i^{ -1 }\,
 h\,
 b\,
 z\,
 i^{ -1 }\,
 c^{ -1 }\rightarrow z\,
 i^{ -1 }\,
 c^{ -1 } - b^{ -1 }\,
 a\,
 f + e^{ -1 }\,
 d\,
 f + e^{ -1 }\,
 d\,
 g\,
 i^{ -1 }\,
 c^{ -1 } + b^{ -1 }\,
 a\,
 g\,
 d\,
 g\,
 i^{ -1 }\,
 c^{ -1 } + b^{ -1 }\,
 a\,
 g\,
 e\,
 z\,
 i^{ -1 }\,
 c^{ -1 } - e^{ -1 }\,
 d\,
 g\,
 d\,
 g\,
 i^{ -1 }\,
 c^{ -1 } - e^{ -1 }\,
 d\,
 g\,
 e\,
 z\,
 i^{ -1 }\,
 c^{ -1 } + e^{ -1 }\,
 i^{ -1 }\,
 h\,
 a\,
 g\,
 i^{ -1 }\,
 c^{ -1 } + e^{ -1 }\,
 i^{ -1 }\,
 h\,
 b\,
 z\,
 i^{ -1 }\,
 c^{ -1 } + b^{ -1 }\,
 a\,
 g\,
 i^{ -1 }\,
 h\,
 a\,
 g\,
 i^{ -1 }\,
 c^{ -1 } + b^{ -1 }\,
 a\,
 g\,
 i^{ -1 }\,
 h\,
 b\,
 z\,
 i^{ -1 }\,
 c^{ -1 } - e^{ -1 }\,
 d\,
 g\,
 i^{ -1 }\,
 h\,
 a\,
 g\,
 i^{ -1 }\,
 c^{ -1 }
$
\end{minipage}\medskip \\
\begin{minipage}{6in}
$
i^{ -1 }\,
 h\,
 b\,
 z\,
 i^{ -1 }\,
 c^{ -1 }\,
 k^{ -1 }\,
 z\,
 d\rightarrow d + i^{ -1 }\,
 h\,
 a - d\,
 g\,
 i^{ -1 }\,
 c^{ -1 }\,
 k^{ -1 }\,
 j\,
 a - d\,
 g\,
 i^{ -1 }\,
 c^{ -1 }\,
 k^{ -1 }\,
 z\,
 d - e\,
 z\,
 i^{ -1 }\,
 c^{ -1 }\,
 k^{ -1 }\,
 j\,
 a - e\,
 z\,
 i^{ -1 }\,
 c^{ -1 }\,
 k^{ -1 }\,
 z\,
 d - i^{ -1 }\,
 h\,
 a\,
 g\,
 i^{ -1 }\,
 c^{ -1 }\,
 k^{ -1 }\,
 j\,
 a - i^{ -1 }\,
 h\,
 a\,
 g\,
 i^{ -1 }\,
 c^{ -1 }\,
 k^{ -1 }\,
 z\,
 d - i^{ -1 }\,
 h\,
 b\,
 z\,
 i^{ -1 }\,
 c^{ -1 }\,
 k^{ -1 }\,
 j\,
 a
$
\end{minipage}\medskip \\
\begin{minipage}{6in}
$
i^{ -1 }\,
 h\,
 b\,
 z\,
 i^{ -1 }\,
 c^{ -1 }\,
 k^{ -1 }\,
 z\,
 e\rightarrow e + i^{ -1 }\,
 h\,
 b + d\,
 g\,
 i^{ -1 }\,
 c^{ -1 }\,
 k^{ -1 } + e\,
 z\,
 i^{ -1 }\,
 c^{ -1 }\,
 k^{ -1 } - d\,
 g\,
 i^{ -1 }\,
 c^{ -1 }\,
 k^{ -1 }\,
 j\,
 b - d\,
 g\,
 i^{ -1 }\,
 c^{ -1 }\,
 k^{ -1 }\,
 z\,
 e - e\,
 z\,
 i^{ -1 }\,
 c^{ -1 }\,
 k^{ -1 }\,
 j\,
 b - e\,
 z\,
 i^{ -1 }\,
 c^{ -1 }\,
 k^{ -1 }\,
 z\,
 e + i^{ -1 }\,
 h\,
 a\,
 g\,
 i^{ -1 }\,
 c^{ -1 }\,
 k^{ -1 } + i^{ -1 }\,
 h\,
 b\,
 z\,
 i^{ -1 }\,
 c^{ -1 }\,
 k^{ -1 } - i^{ -1 }\,
 h\,
 a\,
 g\,
 i^{ -1 }\,
 c^{ -1 }\,
 k^{ -1 }\,
 j\,
 b - i^{ -1 }\,
 h\,
 a\,
 g\,
 i^{ -1 }\,
 c^{ -1 }\,
 k^{ -1 }\,
 z\,
 e - i^{ -1 }\,
 h\,
 b\,
 z\,
 i^{ -1 }\,
 c^{ -1 }\,
 k^{ -1 }\,
 j\,
 b
$
\end{minipage}\medskip \\
\begin{minipage}{6in}
$
k^{ -1 }\,
 z\,
 d\,
 f\,
 k^{ -1 }\,
 j\,
 a\,
 d^{ -1 }\,
 e\rightarrow -1\,
 c\,
 h\,
 b - k^{ -1 }\,
 j\,
 b - k^{ -1 }\,
 z\,
 e + c\,
 h\,
 a\,
 d^{ -1 }\,
 e - k^{ -1 }\,
 j\,
 a\,
 f\,
 k^{ -1 } + k^{ -1 }\,
 j\,
 b\,
 j\,
 b - k^{ -1 }\,
 z\,
 d\,
 f\,
 k^{ -1 } + k^{ -1 }\,
 z\,
 e\,
 j\,
 b + k^{ -1 }\,
 j\,
 a\,
 f\,
 k^{ -1 }\,
 j\,
 b - k^{ -1 }\,
 j\,
 b\,
 j\,
 a\,
 d^{ -1 }\,
 e + k^{ -1 }\,
 z\,
 d\,
 f\,
 k^{ -1 }\,
 j\,
 b - k^{ -1 }\,
 z\,
 e\,
 j\,
 a\,
 d^{ -1 }\,
 e - k^{ -1 }\,
 j\,
 a\,
 f\,
 k^{ -1 }\,
 j\,
 a\,
 d^{ -1 }\,
 e
$
\end{minipage}\medskip \\
\begin{minipage}{6in}
$
k^{ -1 }\,
 z\,
 e\,
 j\,
 a\,
 f\,
 k^{ -1 }\,
 z\,
 d\rightarrow -1\,
 c\,
 h\,
 a\,
 g\,
 d + c\,
 h\,
 b\,
 j\,
 a + k^{ -1 }\,
 j\,
 b\,
 j\,
 a + k^{ -1 }\,
 z\,
 d\,
 g\,
 d + c\,
 h\,
 a\,
 f\,
 k^{ -1 }\,
 j\,
 a + c\,
 h\,
 a\,
 f\,
 k^{ -1 }\,
 z\,
 d - k^{ -1 }\,
 j\,
 a\,
 f\,
 c\,
 h\,
 a + k^{ -1 }\,
 j\,
 a\,
 f\,
 k^{ -1 }\,
 j\,
 a + k^{ -1 }\,
 j\,
 a\,
 f\,
 k^{ -1 }\,
 z\,
 d - k^{ -1 }\,
 j\,
 a\,
 g\,
 d\,
 g\,
 d + k^{ -1 }\,
 j\,
 a\,
 g\,
 e\,
 j\,
 a + k^{ -1 }\,
 j\,
 b\,
 j\,
 a\,
 g\,
 d - k^{ -1 }\,
 j\,
 b\,
 j\,
 b\,
 j\,
 a - k^{ -1 }\,
 z\,
 d\,
 f\,
 c\,
 h\,
 a - k^{ -1 }\,
 z\,
 d\,
 g\,
 d\,
 g\,
 d + k^{ -1 }\,
 z\,
 d\,
 g\,
 e\,
 j\,
 a + k^{ -1 }\,
 z\,
 e\,
 j\,
 a\,
 g\,
 d - k^{ -1 }\,
 z\,
 e\,
 j\,
 b\,
 j\,
 a + k^{ -1 }\,
 j\,
 a\,
 g\,
 d\,
 f\,
 k^{ -1 }\,
 j\,
 a + k^{ -1 }\,
 j\,
 a\,
 g\,
 d\,
 f\,
 k^{ -1 }\,
 z\,
 d - k^{ -1 }\,
 j\,
 b\,
 j\,
 a\,
 f\,
 k^{ -1 }\,
 j\,
 a - k^{ -1 }\,
 j\,
 b\,
 j\,
 a\,
 f\,
 k^{ -1 }\,
 z\,
 d + k^{ -1 }\,
 z\,
 d\,
 g\,
 d\,
 f\,
 k^{ -1 }\,
 j\,
 a + k^{ -1 }\,
 z\,
 d\,
 g\,
 d\,
 f\,
 k^{ -1 }\,
 z\,
 d - k^{ -1 }\,
 z\,
 e\,
 j\,
 a\,
 f\,
 k^{ -1 }\,
 j\,
 a
$
\end{minipage}\medskip \\
\begin{minipage}{6in}
$
k^{ -1 }\,
 z\,
 e\,
 j\,
 a\,
 f\,
 k^{ -1 }\,
 z\,
 e\rightarrow -1\,
 c\,
 h\,
 b - k^{ -1 }\,
 j\,
 b - c\,
 h\,
 a\,
 f\,
 k^{ -1 } - c\,
 h\,
 a\,
 g\,
 e + c\,
 h\,
 b\,
 j\,
 b - k^{ -1 }\,
 j\,
 a\,
 f\,
 k^{ -1 } - k^{ -1 }\,
 j\,
 a\,
 g\,
 e + 2\,
 k^{ -1 }\,
 j\,
 b\,
 j\,
 b + k^{ -1 }\,
 z\,
 e\,
 j\,
 b + c\,
 h\,
 a\,
 f\,
 k^{ -1 }\,
 j\,
 b + c\,
 h\,
 a\,
 f\,
 k^{ -1 }\,
 z\,
 e - k^{ -1 }\,
 j\,
 a\,
 f\,
 c\,
 h\,
 b + k^{ -1 }\,
 j\,
 a\,
 f\,
 k^{ -1 }\,
 j\,
 b + k^{ -1 }\,
 j\,
 a\,
 f\,
 k^{ -1 }\,
 z\,
 e - k^{ -1 }\,
 j\,
 a\,
 g\,
 d\,
 f\,
 k^{ -1 } \\
- 
k^{ -1 }\,
 j\,
 a\,
 g\,
 d\,
 g\,
 e + k^{ -1 }\,
 j\,
 a\,
 g\,
 e\,
 j\,
 b + k^{ -1 }\,
 j\,
 b\,
 j\,
 a\,
 f\,
 k^{ -1 } + k^{ -1 }\,
 j\,
 b\,
 j\,
 a\,
 g\,
 e - k^{ -1 }\,
 j\,
 b\,
 j\,
 b\,
 j\,
 b - k^{ -1 }\,
 z\,
 d\,
 f\,
 c\,
 h\,
 b - k^{ -1 }\,
 z\,
 d\,
 g\,
 d\,
 f\,
 k^{ -1 } - k^{ -1 }\,
 z\,
 d\,
 g\,
 d\,
 g\,
 e + k^{ -1 }\,
 z\,
 d\,
 g\,
 e\,
 j\,
 b + k^{ -1 }\,
 z\,
 e\,
 j\,
 a\,
 f\,
 k^{ -1 } + k^{ -1 }\,
 z\,
 e\,
 j\,
 a\,
 g\,
 e - k^{ -1 }\,
 z\,
 e\,
 j\,
 b\,
 j\,
 b + k^{ -1 }\,
 j\,
 a\,
 g\,
 d\,
 f\,
 k^{ -1 }\,
 j\,
 b + k^{ -1 }\,
 j\,
 a\,
 g\,
 d\,
 f\,
 k^{ -1 }\,
 z\,
 e - k^{ -1 }\,
 j\,
 b\,
 j\,
 a\,
 f\,
 k^{ -1 }\,
 j\,
 b - k^{ -1 }\,
 j\,
 b\,
 j\,
 a\,
 f\,
 k^{ -1 }\,
 z\,
 e + k^{ -1 }\,
 z\,
 d\,
 g\,
 d\,
 f\,
 k^{ -1 }\,
 j\,
 b + k^{ -1 }\,
 z\,
 d\,
 g\,
 d\,
 f\,
 k^{ -1 }\,
 z\,
 e - k^{ -1 }\,
 z\,
 e\,
 j\,
 a\,
 f\,
 k^{ -1 }\,
 j\,
 b
$
\end{minipage}\medskip \\
\begin{minipage}{6in}
$
k^{ -1 }\,
 z\,
 e\,
 j\,
 b\,
 j\,
 b\,
 z\,
 i^{ -1 }\rightarrow c\,
 h\,
 a\,
 f\,
 c + c\,
 h\,
 a\,
 g\,
 i^{ -1 } - k^{ -1 }\,
 z\,
 d\,
 f\,
 c - k^{ -1 }\,
 z\,
 d\,
 g\,
 i^{ -1 } - c\,
 h\,
 a\,
 g\,
 d\,
 g\,
 i^{ -1 } - c\,
 h\,
 a\,
 g\,
 e\,
 z\,
 i^{ -1 } + c\,
 h\,
 b\,
 j\,
 a\,
 g\,
 i^{ -1 } + c\,
 h\,
 b\,
 j\,
 b\,
 z\,
 i^{ -1 } + k^{ -1 }\,
 j\,
 a\,
 g\,
 d\,
 f\,
 c + k^{ -1 }\,
 j\,
 a\,
 g\,
 d\,
 g\,
 i^{ -1 } - k^{ -1 }\,
 j\,
 b\,
 j\,
 a\,
 f\,
 c + k^{ -1 }\,
 j\,
 b\,
 j\,
 b\,
 z\,
 i^{ -1 } + k^{ -1 }\,
 z\,
 d\,
 g\,
 d\,
 f\,
 c + 2\,
 k^{ -1 }\,
 z\,
 d\,
 g\,
 d\,
 g\,
 i^{ -1 } + k^{ -1 }\,
 z\,
 d\,
 g\,
 e\,
 z\,
 i^{ -1 } - k^{ -1 }\,
 z\,
 e\,
 j\,
 a\,
 f\,
 c - k^{ -1 }\,
 z\,
 e\,
 j\,
 a\,
 g\,
 i^{ -1 } - k^{ -1 }\,
 j\,
 a\,
 f\,
 c\,
 h\,
 a\,
 g\,
 i^{ -1 } - k^{ -1 }\,
 j\,
 a\,
 f\,
 c\,
 h\,
 b\,
 z\,
 i^{ -1 } - k^{ -1 }\,
 j\,
 a\,
 g\,
 d\,
 g\,
 d\,
 g\,
 i^{ -1 } - k^{ -1 }\,
 j\,
 a\,
 g\,
 d\,
 g\,
 e\,
 z\,
 i^{ -1 } + k^{ -1 }\,
 j\,
 a\,
 g\,
 e\,
 j\,
 a\,
 g\,
 i^{ -1 } + k^{ -1 }\,
 j\,
 a\,
 g\,
 e\,
 j\,
 b\,
 z\,
 i^{ -1 } + k^{ -1 }\,
 j\,
 b\,
 j\,
 a\,
 g\,
 d\,
 g\,
 i^{ -1 } + k^{ -1 }\,
 j\,
 b\,
 j\,
 a\,
 g\,
 e\,
 z\,
 i^{ -1 } - k^{ -1 }\,
 j\,
 b\,
 j\,
 b\,
 j\,
 a\,
 g\,
 i^{ -1 } - k^{ -1 }\,
 j\,
 b\,
 j\,
 b\,
 j\,
 b\,
 z\,
 i^{ -1 } - k^{ -1 }\,
 z\,
 d\,
 f\,
 c\,
 h\,
 a\,
 g\,
 i^{ -1 } - k^{ -1 }\,
 z\,
 d\,
 f\,
 c\,
 h\,
 b\,
 z\,
 i^{ -1 } - k^{ -1 }\,
 z\,
 d\,
 g\,
 d\,
 g\,
 d\,
 g\,
 i^{ -1 } - k^{ -1 }\,
 z\,
 d\,
 g\,
 d\,
 g\,
 e\,
 z\,
 i^{ -1 } + k^{ -1 }\,
 z\,
 d\,
 g\,
 e\,
 j\,
 a\,
 g\,
 i^{ -1 } + k^{ -1 }\,
 z\,
 d\,
 g\,
 e\,
 j\,
 b\,
 z\,
 i^{ -1 } + k^{ -1 }\,
 z\,
 e\,
 j\,
 a\,
 g\,
 d\,
 g\,
 i^{ -1 } + k^{ -1 }\,
 z\,
 e\,
 j\,
 a\,
 g\,
 e\,
 z\,
 i^{ -1 } - k^{ -1 }\,
 z\,
 e\,
 j\,
 b\,
 j\,
 a\,
 g\,
 i^{ -1 }
$
\end{minipage}\\
\vspace{10pt}
\label{end-digested}

\noindent
The time for preliminaries was 0:00:01\\
The time for NCMakeGB 1 was 0:00:00\\
The time for Remove Redundant 1 was 0:00:00\\
The time for the main NCMakeGB was 0:00:05\\
The time for Remove Redundant 2 was 0:00:00\\
The time for reducing unknowns was 0:00:01\\
The time for clean up basis was 0:00:02\\
The time for SmallBasis was 0:00:01\\
The time for CreateCategories was 0:01:07\\
The time for NCCV was 0:00:00\\
The time for RegularOutput was 0:00:38\\
The time for everything so far was 0:01:58\\


\def\CellGroup{\bgroup}
\def\endCellGroup{\egroup}



\section{ Appendix 3 - Find A Smaller Basis   }
\label{sec-attach2}

\mathin
Input: 

SetNonCommutative[a,b,c,d,e,f,g,h,i,w,x,y,z,j,u,t,k,v,\ 
Inv[a],Inv[b],Inv[c],Inv[d],Inv[e],Inv[f],Inv[g],Inv[h],
Inv[i],Inv[j], Inv[k]  ];

(*Here are the relations which we like and would like to retain*)
hopepolys =   \{-1 + a**Inv[a], -1 + b**Inv[b], -1 + c**Inv[c], 
   -1 + d**Inv[d], -1 + e**Inv[e], -1 + f**Inv[f], 
   -1 + g**Inv[g], -1 + h**Inv[h], -1 + i**Inv[i], 
   -1 + j**Inv[j], 3-1 + k**Inv[k], -1 + Inv[a]**a, \
   -1 + Inv[b]**b, -1 + Inv[c]**c, -1 + Inv[d]**d, -1 + Inv[e]**e, \
   -1 + Inv[f]**f, -1 + Inv[g]**g, -1 + Inv[h]**h, -1 + Inv[i]**i, 
   -1 + Inv[j]**j, -1 + Inv[k]**k, t + a**g**Inv[i] + b**z**Inv[i], 
   w - Inv[i] + d**g**Inv[i] + e**z**Inv[i], 
   -z + j**a**g + j**b**z - k**c**i + z**d**g + z**e**z, 
   u + Inv[k]**j**a + Inv[k]**z**d, v - Inv[k] + Inv[k]**j**b + Inv[k]**z**e, 
   x - Inv[a] - f**Inv[k]**j + g**d**Inv[a] - f**Inv[k]**z**d**Inv[a], 
   y - Inv[c]**Inv[k]**j**a**f - Inv[c]**Inv[k]**j**b**k - 
    Inv[c]**Inv[k]**z**d**f - Inv[c]**Inv[k]**z**e**k  
\}

(* Here are the relations resulting from the original "discovery" run *)
  allpolys = \{-1 + a**Inv[a], -1 + b**Inv[b], -1 + c**Inv[c], -1 + d**Inv[d], 
   -1 + e**Inv[e], -1 + f**Inv[f], -1 + g**Inv[g], -1 + h**Inv[h], 
   -1 + i**Inv[i], -1 + j**Inv[j], -1 + k**Inv[k], -1 + Inv[a]**a, 
   -1 + Inv[b]**b, -1 + Inv[c]**c, -1 + Inv[d]**d, -1 + Inv[e]**e, 
   -1 + Inv[f]**f, -1 + Inv[g]**g, -1 + Inv[h]**h, -1 + Inv[i]**i, 
   -1 + Inv[j]**j, -1 + Inv[k]**k, t + a**g**Inv[i] + b**z**Inv[i], 
   w - Inv[i] + d**g**Inv[i] + e**z**Inv[i], 
   -z + j**a**g + j**b**z - k**c**i + z**d**g + z**e**z, 
   u + Inv[k]**j**a + Inv[k]**z**d, v - Inv[k] + Inv[k]**j**b + Inv[k]**z**e, 
   -Inv[a]**Inv[h] - Inv[d]**Inv[i] + Inv[a]**b**j**Inv[h] - 
   Inv[a]**b**z**Inv[i] - Inv[d]**e**j**Inv[h] + Inv[d]**e**z**Inv[i], 
   -Inv[f]**Inv[a] - Inv[k]**Inv[b] + Inv[f]**g**d**Inv[a] - 
   Inv[f]**g**e**Inv[b] - Inv[k]**z**d**Inv[a] + Inv[k]**z**e**Inv[b], 
   x - Inv[a] - f**Inv[k]**j + g**d**Inv[a] - f**Inv[k]**z**d**Inv[a], 
   y - Inv[c]**Inv[k]**j**a**f - Inv[c]**Inv[k]**j**b**k - 
   Inv[c]**Inv[k]**z**d**f - Inv[c]**Inv[k]**z**e**k, 
   Inv[k] + Inv[f]**g**e + Inv[f]**Inv[a]**b - Inv[k]**z**e - 
   Inv[f]**g**d**Inv[a]**b + Inv[k]**z**d**Inv[a]**b, 
   h**a**g**Inv[i]**Inv[c] + h**b**z**Inv[i]**Inv[c] + 
   i**d**g**Inv[i]**Inv[c] + i**e**z**Inv[i]**Inv[c] - 
   Inv[c]**Inv[k]**j**a**f - Inv[c]**Inv[k]**j**b**k - 
   Inv[c]**Inv[k]**z**d**f - Inv[c]**Inv[k]**z**e**k, 
\endmathin

\mathin 
   -c**h**a**g**Inv[i] - c**h**b**z**Inv[i] - c**i**d**g**Inv[i] - 
   c**i**e**z**Inv[i] + Inv[k]**j**a**f**c + Inv[k]**j**b**k**c + 
   Inv[k]**z**d**f**c + Inv[k]**z**e**k**c, 
   b**j + a**f**Inv[k]**j - a**g**d**Inv[a] - a**g**Inv[i]**h - 
   b**z**Inv[i]**h + a**f**Inv[k]**z**d**Inv[a], 
   d**Inv[a] + e**j + Inv[i]**h + d**f**Inv[k]**j - d**g**d**Inv[a] - 
   d**g**Inv[i]**h - e**z**Inv[i]**h + d**f**Inv[k]**z**d**Inv[a], 
   -g**Inv[i] + f**Inv[k]**j**Inv[h] - g**d**Inv[a]**Inv[h] + 
   Inv[a]**b**j**Inv[h] - Inv[a]**b**z**Inv[i] + 
   f**Inv[k]**z**d**Inv[a]**Inv[h], 
   -h**a - i**d + h**a**g**d - h**b**j**a + i**d**g**d - i**e**j**a - 
   h**a**f**Inv[k]**j**a - h**a**f**Inv[k]**z**d + h**a**g**Inv[i]**h**a + 
   h**b**z**Inv[i]**h**a - i**d**f**Inv[k]**j**a - i**d**f**Inv[k]**z**d + 
   i**d**g**Inv[i]**h**a + i**e**z**Inv[i]**h**a, 
   h**a**f**Inv[k] + h**a**g**e - h**b**j**b + i**d**f**Inv[k] + 
   i**d**g**e - i**e**j**b - h**a**f**Inv[k]**j**b - 
   h**a**f**Inv[k]**z**e + h**a**g**Inv[i]**h**b + h**b**z**Inv[i]**h**b - 
   i**d**f**Inv[k]**j**b - i**d**f**Inv[k]**z**e + i**d**g**Inv[i]**h**b + 
   i**e**z**Inv[i]**h**b, -j**a**f**c - j**b**k**c - z**d**f**c - 
   z**e**k**c + k**c**h**a**g**Inv[i] + k**c**h**b**z**Inv[i] + 
   k**c**i**d**g**Inv[i] + k**c**i**e**z**Inv[i], 
   z**d - j**a**g**d + j**b**j**a - k**c**h**a - z**d**g**d + z**e**j**a + 
   j**a**f**Inv[k]**j**a + j**a**f**Inv[k]**z**d + z**d**f**Inv[k]**j**a + 
   z**d**f**Inv[k]**z**d, -j**b - j**a**f**Inv[k] - j**a**g**e + 
   j**b**j**b - k**c**h**b - z**d**f**Inv[k] - z**d**g**e + z**e**j**b + 
   j**a**f**Inv[k]**j**b + j**a**f**Inv[k]**z**e + z**d**f**Inv[k]**j**b + 
   z**d**f**Inv[k]**z**e, -a**f**c + b**z**Inv[i] + a**g**d**g**Inv[i] + 
   a**g**e**z**Inv[i] - b**j**a**g**Inv[i] - b**j**b**z**Inv[i] + 
   a**g**Inv[i]**h**a**g**Inv[i] + a**g**Inv[i]**h**b**z**Inv[i] + 
   b**z**Inv[i]**h**a**g**Inv[i] + b**z**Inv[i]**h**b**z**Inv[i], 
   -d**f**c - d**g**Inv[i] + d**g**d**g**Inv[i] + d**g**e**z**Inv[i] - 
   e**j**a**g**Inv[i] - e**j**b**z**Inv[i] - Inv[i]**h**a**g**Inv[i] - 
   Inv[i]**h**b**z**Inv[i] + d**g**Inv[i]**h**a**g**Inv[i] + 
   d**g**Inv[i]**h**b**z**Inv[i] + e**z**Inv[i]**h**a**g**Inv[i] + 
   e**z**Inv[i]**h**b**z**Inv[i], 
   Inv[k]**z - c**h**a**Inv[d] - Inv[k]**j**a**g - Inv[k]**z**d**g + 
   Inv[k]**j**a**f**Inv[k]**z + Inv[k]**j**b**j**a**Inv[d] + 
   Inv[k]**z**d**f**Inv[k]**z + Inv[k]**z**e**j**a**Inv[d] + 
   Inv[k]**j**a**f**Inv[k]**j**a**Inv[d] + 
   Inv[k]**z**d**f**Inv[k]**j**a**Inv[d], 
   -Inv[k]**z + c**h**a**Inv[d] - c**h**b**Inv[e] - Inv[k]**j**b**Inv[e] - 
   Inv[k]**j**a**f**Inv[k]**Inv[e] - Inv[k]**j**b**j**a**Inv[d] + 
   Inv[k]**j**b**j**b**Inv[e] - Inv[k]**z**d**f**Inv[k]**Inv[e] - 
   Inv[k]**z**e**j**a**Inv[d] + Inv[k]**z**e**j**b**Inv[e] - 
   Inv[k]**j**a**f**Inv[k]**j**a**Inv[d] +
   Inv[k]**j**a**f**Inv[k]**j**b**Inv[e] 
   - Inv[k]**z**d**f**Inv[k]**j**a**Inv[d] + 
   Inv[k]**z**d**f**Inv[k]**j**b**Inv[e], 
   z**Inv[i]**Inv[c] - Inv[b]**a**f - j**a**g**Inv[i]**Inv[c] - 
\endmathin
\mathin
   j**b**z**Inv[i]**Inv[c] + z**Inv[i]**h**a**g**Inv[i]**Inv[c] + 
   z**Inv[i]**h**b**z**Inv[i]**Inv[c] + Inv[b]**a**g**d**g**Inv[i]**Inv[c] + 
   Inv[b]**a**g**e**z**Inv[i]**Inv[c] + 
   Inv[b]**a**g**Inv[i]**h**a**g**Inv[i]**Inv[c] + 
   Inv[b]**a**g**Inv[i]**h**b**z**Inv[i]**Inv[c], 
   -z**Inv[i]**Inv[c] + Inv[b]**a**f - Inv[e]**d**f - 
   Inv[e]**d**g**Inv[i]**Inv[c] - Inv[b]**a**g**d**g**Inv[i]**Inv[c] - 
   Inv[b]**a**g**e**z**Inv[i]**Inv[c] + Inv[e]**d**g**d**g**Inv[i]**Inv[c] + 
   Inv[e]**d**g**e**z**Inv[i]**Inv[c] - 
   Inv[e]**Inv[i]**h**a**g**Inv[i]**Inv[c] - 
   Inv[e]**Inv[i]**h**b**z**Inv[i]**Inv[c] - 
   Inv[b]**a**g**Inv[i]**h**a**g**Inv[i]**Inv[c] - 
   Inv[b]**a**g**Inv[i]**h**b**z**Inv[i]**Inv[c] + 
   Inv[e]**d**g**Inv[i]**h**a**g**Inv[i]**Inv[c] + 
   Inv[e]**d**g**Inv[i]**h**b**z**Inv[i]**Inv[c], 
   -d - Inv[i]**h**a + d**g**Inv[i]**Inv[c]**Inv[k]**j**a + 
   d**g**Inv[i]**Inv[c]**Inv[k]**z**d + e**z**Inv[i]**Inv[c]**Inv[k]**j**a + 
   e**z**Inv[i]**Inv[c]**Inv[k]**z**d + 
   Inv[i]**h**a**g**Inv[i]**Inv[c]**Inv[k]**j**a + 
   Inv[i]**h**a**g**Inv[i]**Inv[c]**Inv[k]**z**d + 
   Inv[i]**h**b**z**Inv[i]**Inv[c]**Inv[k]**j**a + 
   Inv[i]**h**b**z**Inv[i]**Inv[c]**Inv[k]**z**d, 
   -e - Inv[i]**h**b - d**g**Inv[i]**Inv[c]**Inv[k] - 
   e**z**Inv[i]**Inv[c]**Inv[k] + d**g**Inv[i]**Inv[c]**Inv[k]**j**b + 
   d**g**Inv[i]**Inv[c]**Inv[k]**z**e + e**z**Inv[i]**Inv[c]**Inv[k]**j**b + 
   e**z**Inv[i]**Inv[c]**Inv[k]**z**e - 
   Inv[i]**h**a**g**Inv[i]**Inv[c]**Inv[k] - 
   Inv[i]**h**b**z**Inv[i]**Inv[c]**Inv[k] + 
   Inv[i]**h**a**g**Inv[i]**Inv[c]**Inv[k]**j**b + 
   Inv[i]**h**a**g**Inv[i]**Inv[c]**Inv[k]**z**e + 
   Inv[i]**h**b**z**Inv[i]**Inv[c]**Inv[k]**j**b + 
   Inv[i]**h**b**z**Inv[i]**Inv[c]**Inv[k]**z**e, 
   c**h**b + Inv[k]**j**b + Inv[k]**z**e - c**h**a**Inv[d]**e + 
   Inv[k]**j**a**f**Inv[k] - Inv[k]**j**b**j**b + Inv[k]**z**d**f**Inv[k] - 
   Inv[k]**z**e**j**b - Inv[k]**j**a**f**Inv[k]**j**b + 
   Inv[k]**j**b**j**a**Inv[d]**e - Inv[k]**z**d**f**Inv[k]**j**b + 
   Inv[k]**z**e**j**a**Inv[d]**e + Inv[k]**j**a**f**Inv[k]**j**a**Inv[d]**e + 
   Inv[k]**z**d**f**Inv[k]**j**a**Inv[d]**e, 
   c**h**a**g**d - c**h**b**j**a - Inv[k]**j**b**j**a - Inv[k]**z**d**g**d - 
   c**h**a**f**Inv[k]**j**a - c**h**a**f**Inv[k]**z**d + 
   Inv[k]**j**a**f**c**h**a - Inv[k]**j**a**f**Inv[k]**j**a - 
   Inv[k]**j**a**f**Inv[k]**z**d + Inv[k]**j**a**g**d**g**d - 
   Inv[k]**j**a**g**e**j**a - Inv[k]**j**b**j**a**g**d + 
   Inv[k]**j**b**j**b**j**a + Inv[k]**z**d**f**c**h**a + 
   Inv[k]**z**d**g**d**g**d - Inv[k]**z**d**g**e**j**a - 
   Inv[k]**z**e**j**a**g**d + Inv[k]**z**e**j**b**j**a - 
   Inv[k]**j**a**g**d**f**Inv[k]**j**a - Inv[k]**j**a**g**d**f**Inv[k]**z**d + 
   Inv[k]**j**b**j**a**f**Inv[k]**j**a + Inv[k]**j**b**j**a**f**Inv[k]**z**d - 
   Inv[k]**z**d**g**d**f**Inv[k]**j**a - Inv[k]**z**d**g**d**f**Inv[k]**z**d + 
   Inv[k]**z**e**j**a**f**Inv[k]**j**a + Inv[k]**z**e**j**a**f**Inv[k]**z**d,
\endmathin
\mathin
   c**h**b + Inv[k]**j**b + c**h**a**f**Inv[k] + c**h**a**g**e - 
   c**h**b**j**b + Inv[k]**j**a**f**Inv[k] + Inv[k]**j**a**g**e - 
   2*Inv[k]**j**b**j**b - Inv[k]**z**e**j**b - c**h**a**f**Inv[k]**j**b - 
   c**h**a**f**Inv[k]**z**e + Inv[k]**j**a**f**c**h**b - 
   Inv[k]**j**a**f**Inv[k]**j**b - Inv[k]**j**a**f**Inv[k]**z**e + 
   Inv[k]**j**a**g**d**f**Inv[k] + Inv[k]**j**a**g**d**g**e - 
   Inv[k]**j**a**g**e**j**b - Inv[k]**j**b**j**a**f**Inv[k] - 
   Inv[k]**j**b**j**a**g**e + Inv[k]**j**b**j**b**j**b + 
   Inv[k]**z**d**f**c**h**b + Inv[k]**z**d**g**d**f**Inv[k] + 
   Inv[k]**z**d**g**d**g**e - Inv[k]**z**d**g**e**j**b - 
   Inv[k]**z**e**j**a**f**Inv[k] - Inv[k]**z**e**j**a**g**e + 
   Inv[k]**z**e**j**b**j**b - Inv[k]**j**a**g**d**f**Inv[k]**j**b - 
   Inv[k]**j**a**g**d**f**Inv[k]**z**e + Inv[k]**j**b**j**a**f**Inv[k]**j**b + 
   Inv[k]**j**b**j**a**f**Inv[k]**z**e - Inv[k]**z**d**g**d**f**Inv[k]**j**b - 
   Inv[k]**z**d**g**d**f**Inv[k]**z**e + 
   Inv[k]**z**e**j**a**f**Inv[k]**j**b + Inv[k]**z**e**j**a**f**Inv[k]**z**e,
   -c**h**a**f**c - c**h**a**g**Inv[i] + Inv[k]**z**d**f**c + 
   Inv[k]**z**d**g**Inv[i] + c**h**a**g**d**g**Inv[i] + 
   c**h**a**g**e**z**Inv[i] - c**h**b**j**a**g**Inv[i] - 
   c**h**b**j**b**z**Inv[i] - Inv[k]**j**a**g**d**f**c - 
   Inv[k]**j**a**g**d**g**Inv[i] + Inv[k]**j**b**j**a**f**c - 
   Inv[k]**j**b**j**b**z**Inv[i] - Inv[k]**z**d**g**d**f**c - 
   2*Inv[k]**z**d**g**d**g**Inv[i] - Inv[k]**z**d**g**e**z**Inv[i] + 
   Inv[k]**z**e**j**a**f**c + Inv[k]**z**e**j**a**g**Inv[i] + 
   Inv[k]**j**a**f**c**h**a**g**Inv[i] + Inv[k]**j**a**f**c**h**b**z**Inv[i] + 
   Inv[k]**j**a**g**d**g**d**g**Inv[i] + Inv[k]**j**a**g**d**g**e**z**Inv[i] - 
   Inv[k]**j**a**g**e**j**a**g**Inv[i] - Inv[k]**j**a**g**e**j**b**z**Inv[i] - 
   Inv[k]**j**b**j**a**g**d**g**Inv[i] - Inv[k]**j**b**j**a**g**e**z**Inv[i] + 
   Inv[k]**j**b**j**b**j**a**g**Inv[i] + Inv[k]**j**b**j**b**j**b**z**Inv[i] + 
   Inv[k]**z**d**f**c**h**a**g**Inv[i] + Inv[k]**z**d**f**c**h**b**z**Inv[i] + 
   Inv[k]**z**d**g**d**g**d**g**Inv[i] + Inv[k]**z**d**g**d**g**e**z**Inv[i] - 
   Inv[k]**z**d**g**e**j**a**g**Inv[i] - Inv[k]**z**d**g**e**j**b**z**Inv[i] - 
   Inv[k]**z**e**j**a**g**d**g**Inv[i] - Inv[k]**z**e**j**a**g**e**z**Inv[i] + 
   Inv[k]**z**e**j**b**j**a**g**Inv[i] + Inv[k]**z**e**j**b**j**b**z**Inv[i]\}

(* Here we set the order to be strictly graded lex *) 

SetKnowns[a,Inv[a],b,Inv[b],c,Inv[c],d,Inv[d],e, Inv[e], f, Inv[f], g, Inv[g], h, Inv[h],i, Inv[i],j, Inv[j], k, Inv[k] , z, y,x,t,u,v,w];

SetUnknowns[\{\}];

(* Ask for a small basis retaining the relations we like *) 
SmallBasis[ allpolys, hopepolys, 3 ]
\endmathin

\noindent
Output: 
\label{smallb}
\begin{CellGroup}

\dispSFoutmath{\MathBegin{MathArray}{l}
\{-\Muserfunction{Inv}[a]**\Muserfunction{Inv}[h]-\Muserfunction{Inv}[d]**\Muserfunction{Inv}[i]+  
\Muserfunction{Inv}[a]**b**j**\Muserfunction{Inv}[h]-  \\
\noalign{\vspace{0.5ex}}
\hspace{2.em} \Muserfunction{Inv}[a]**b**z**\Muserfunction{Inv}[i]-  
\Muserfunction{Inv}[d]**e**j**\Muserfunction{Inv}[h]+  \\
\noalign{\vspace{0.5ex}}
\hspace{2.em} \Muserfunction{Inv}[d]**e**z**\Muserfunction{Inv}[i],  \\
\noalign{\vspace{0.5ex}}
\hspace{1.em} -\Muserfunction{Inv}[f]**\Muserfunction{Inv}[a]-\Muserfunction{Inv}[k]**\Muserfunction{Inv}[b]+  \\
\noalign{\vspace{0.5ex}}
\hspace{2.em} \Muserfunction{Inv}[f]**g**d**\Muserfunction{Inv}[a]-  
\Muserfunction{Inv}[f]**g**e**\Muserfunction{Inv}[b]-  \\
\noalign{\vspace{0.5ex}}
\hspace{2.em} \Muserfunction{Inv}[k]**z**d**\Muserfunction{Inv}[a]+  
\Muserfunction{Inv}[k]**z**e**\Muserfunction{Inv}[b]\}\\
\MathEnd{MathArray}}

\end{CellGroup}


\def\CellGroup{\bgroup}
\def\endCellGroup{\egroup}


\section{ Appendix 4 - Confirm Our Relations Imply The Result }
\label{sec-attach3}
\noindent
\mathin
Input:
 
SetNonCommutative[a,b,c,d,e,f,g,h,i,w,x,y,z,j,u,t,k,v,\ 
Inv[a],Inv[b],Inv[c],Inv[d],Inv[e],Inv[f],Inv[g],Inv[h],
Inv[i],Inv[j], Inv[k]  ];

(* Here we create the relations which we wish to imply *)
first = \{\{a,t,b\},\{u,c,v\},\{d,w, e\}\};

second = \{\{ x,f,g\}, \{h,y,i\},\{j,k,z\}\};

oneway =MatMult[first,second] - IdentityMatrix[3];

otherway = 
     MatMult[second,first]-  IdentityMatrix[3]; 

start = Flatten[\{ oneway, otherway  \}];

SetKnowns[a,Inv[a],b,Inv[b],c,Inv[c],d,Inv[d],e, Inv[e], f, Inv[f], g, Inv[g], h, Inv[h],i, Inv[i],j, Inv[j], k, Inv[k] , z, y,x,t,u,v,w];

SetUnknowns[\{\}];

(* hopepolys = *)
(* hopepolys are defined as in Appendix 3 *)

(* Here are the relations found through the Small Basis algorithm *)

newrels = \{-Inv[a]**Inv[h] - Inv[d]**Inv[i] + Inv[a]**b**j**Inv[h] - 
   Inv[a]**b**z**Inv[i] - Inv[d]**e**j**Inv[h] + Inv[d]**e**z**Inv[i], 
  -Inv[f]**Inv[a] - Inv[k]**Inv[b] + Inv[f]**g**d**Inv[a] - 
   Inv[f]**g**e**Inv[b] - Inv[k]**z**d**Inv[a] + Inv[k]**z**e**Inv[b]\}

hopepolys = Join[ hopepolys, newrels ];

(*  We will have our result if the starting relations are elements of
    the ideal generated by the above relations i.e. the GB created with 
    the above relations reduces the starting relations to 0 *)

(*  Note that a 5 iteration GB failed to reduce the starting relations *)
hopeGB = NCMakeGB[ hopepolys , 7 ];

hoperules =  PolyToRule[ hopeGB ];

(* Use the Groebner basis relations to reduce the original relations *)
Reduction[ start , hoperules ]

\endmathin

\begin{CellGroup}
Output:
\dispSFoutmath{\MathBegin{MathArray}{l}
\{0,0,0,0,0,0,0,0,0,0,0,0,0,0,0,0,0,0\}  \\
\MathEnd{MathArray}}

\end{CellGroup}